\documentclass[12pt]{article}
\usepackage[dvips]{graphicx}
\usepackage{amsmath, amssymb}
\usepackage{amsthm}

\newtheorem{theorem}{Theorem}[section]

\newtheorem{proposition}{Proposition}[section]
\newtheorem{definition}{Definition}[section]

\newcommand{\tp}{^{\mathrm t}}

\newlength{\IndentI}
\newlength{\IndentII}
\newlength{\IndentIII}
\newlength{\WidthI}
\newlength{\WidthII}
\newlength{\WidthIII}
\title{Conformal geometry of statistical manifold with application to  sequential estimation}

\author{
Masayuki Kumon\footnote{
masayuki\_kumon@smile.odn.ne.jp}  
\\
Akimichi Takemura
\footnote{
Graduate School of Information Science and Technology,
University of Tokyo, \   
7-3-1 Hongo, Bunkyo-ku, Tokyo 113-8656, JAPAN, \ 
takemura@stat.t.u-tokyo.ac.jp}
\\
Kei Takeuchi
\footnote{
Emeritus, Graduate School of Economics, 
University of Tokyo}
}

\date{February 2011}

\begin{document}
\maketitle

\begin{abstract}
We present a geometrical method for analyzing sequential estimating procedures. 
It is based on the design principle of the second-order efficient sequential estimation provided in Okamoto, Amari and Takeuchi (1991). 
By introducing a dual conformal curvature quantity, we clarify the conditions for the covariance minimization of sequential estimators. 
These conditions are further elabolated for the multidimensional curved exponential family. 
The theoretical results are then numerically examined by using typical statistical models,  
von Mises-Fisher and hyperboloid models. 
\end{abstract}


\section{Introduction}

Sequential estimation continues observations until the observed sample satisfies a certain prescribed criterion. 
Its properties have been shown to be superior on the average to those of nonsequential estimation in which the number of observations is fixed 
a priori. 
Specifically the developments of higher-order asymptotic theory have suggested that the information loss due to the exponential curvature of the statistical model might be recovered by a sequential estimation procedure which makes use of the ancillary statistic. Such an estimator is expected to have a uniformly better characteristic on the average 
(see e.g. S\o rensen (1986)). 
Takeuchi and Akahira (1988) formulated this scheme rigorously and analyzed the higher-order efficiency of sequential estimation procedures in the scalar parameter case (see also Akahira and Takeuchi (1989)). 
They showed that in the sequential case the exponential curvature term in the second-order variance can be eliminated by a second-order efficient estimator, and the maximum likelihood estimator with an appropriate  stopping rule gives such a sequential estimator. 
This also implies that appropriately designed sequential estimators are superior to nonsequential estimators in the asymptotic sense.

Following the work of Takeuchi and Akahira (1988, 1989), Okamoto, Amari and Takeuchi (1991) generalized the results to the multiparameter case by using the geometrical method, and studied characteristics of more general sequential estimation procedures. 
In the nonsequential case, a statistical manifold is uniformly enlarged by $N$ times when we use $N$ observations, keeping the intrinsic features of the manifold unchanged. 
In a sequential estimation procedure with a certain stopping rule, the observed sample size $N$ is a random variable depending on the position of a  statistical manifold. 
This causes a nonuniform expansion of the statistical manifold. 
Such an expansion is called the conformal transformation in geometry, since it changes the scale locally and isotropically but it does not change the shape of a figure (it does not change the orthogonality).  
The result of Takeuchi and Akahira can be interpreted such that it is possible to reduce the exponential curvature of a statistical manifold to zero by a suitable conformal transformation. 
The conformal geometry thus is an adequate framework for the analysis of the sequential inferential procedures if we extend the concept of the conformal transformation to the statistical manifold which is the Riemannian manifold with a dual couple of affine connections.

As a sequel to Okamoto, Amari and Takeuchi (1991), this paper investigates  the sequential estimating procedures from the information geometrical viewpoint. 
The novelty of this paper is the introduction of the dual conformal Weyl-Schouten curvature of a statistical manifold, and this quantity will be proved to play a 
central role when considering the problem of covariance minimization 
under the sequential estimating procedures. 
Information geometry was originated with the work of 
Amari (1985), and it has been establishing a solid status as a mathematical 
methodology for a variety of statistical sciences 
(see e.g. Amari et al (1987), Amari and Nagaoka (2000), 
Kumon (2009, 2010)). 
In line with these developments, the present paper also intends to provide a starting point for studying the conformal geometry of a statistical manifold itself succeeding to the work of Lauritzen (1987).

The organization of the paper is as follows. 
The established known results are cited as propositions, and the results obtained in this paper are stated as theorems. 
In the next section, we prepare some statistical notations and preliminary 
results which will be relevant in this paper. 
In Section 3, we formulate a conformal transformation of a statistical manifold, where a set of dual Weyl-Schouten curvature tensors is  introduced. 
Then we elucidate their implications in the structures of statistical manifolds. In this connection, the meaning of conjugate symmetry is also explained, which is the notion first introduced by Lauritzen (1987).    
In Section 4, the general result in the previous section is used to write out the structure of multidimentional exponential family. 
In Section 5, the general result is used to delineate the structure of multidimensional curved exponential family, where the dual Euler-Schouten curvatures are introduced. Then related notions such as totally exponential umbilic and dual quadric hypersurface are shown to involve key  elements in studying statistical submanifolds.  
In Section 6, the geometrical results obtained in the previous sections are applied to the sequential estimation in a multidimensional curved exponential family, where 
we give a concrete procedure for the covariance minimization. 
In Section 7, the results in Section 6 are numerically examined by two typical curved exponential families called the von Mises-Fisher model and the hyperboloid model. 
Section 8 is devoted to some additional discussions and a perspective of future work.

\section{Preliminaries}

Let us denote by $X(t) = (X_1(t), \dots, X_k(t))^{\tp}$ a $k$-dimentional random process defined on the probability space $[\Omega, {\cal F}, P]$ 
with values in $[E, \cal{E}]$, where $E = \mathbb{R}^k$ and $\cal{E}$ is the $\sigma$-field of all Borel sets in $E$. 
The time parameter $t \in T$ runs over all non-negative integers 
$T = \{0, 1, 2, \dots \}$ or over all non-negative real numbers 
$T =[0, +\infty)$. Moreover, let ${\cal F}_t, t \in T$ denote the $\sigma$-field generated by the random vectors $X(s), s \leqslant t$.

We assume that the probability measure $P$ depends on an unknown parameter $\theta = (\theta^1, \dots, \theta^m)^{\tp} \in \Theta$, $P = P_\theta$, where $\Theta$ is homeomorphic to ${\mathbb R}^m$, 
and we shall consider the case where the following conditions are fullfiled:

(i)\ $X(t)$ is continuous in probability, has stationary  independent increments and $P_\theta(X(0)) = 0) = 1,\ \forall \theta \in \Theta$.

(ii)\ The probability distributions at any time $t$ are dominated by a $\sigma$-finite measure $\mu$ and have the densities 
$f(x, t, \theta)$ with respect to $\mu$
\begin{align}
\label{eq:1}
\frac{dP_\theta}{d\mu}(x, t) = f(x, t, \theta),\quad 
x\in E,\ t \in T,\ \theta \in \Theta.
\end{align}
We say that $\{P_\theta\ |\ \theta\in \Theta\}$ is an $m$-dimensional full regular minimally represented exponential family 
(f.r.m.\ exponential family) when the densities 
(\ref{eq:1}) can be written as
\begin{align}
\label{eq:2}
f_e(x, t, \theta) = \exp \{\theta^ix_i - \psi(\theta)t\},
\end{align}
where $x = (x_1, \dots, x_m)^{\tp}\in {\mathbb R}^m$, $\theta$ is the natural parameter, and 
$\eta = \partial \psi(\theta)/\partial \theta$ is the expectation parameter with $\psi(\theta)$ a smooth (infinitely differentiable) convex function of $\theta$. In the right-hand side of (\ref{eq:2}) and hereafter the Einstein summation convention will be assumed, so that summation will be automatically taken over indices repeated twice in the sense e.g. 
$\theta^ix_i = \sum_{i=1}^m \theta^ix_i$.

Sequential statistical procedures are characterized by a random sample size, where stopping times are used to stop the observations of the process. We denote by $\tau$ an arbitrary stopping time, i.e., a random variable $\tau$ defined on $\Omega$ with values in 
$T\cup \{\infty\}$ and possessing the property $\{\omega \in \Omega: \tau(\omega) \leqslant t\}\in {\cal F}_t,\ \forall t\in T$. 
We consider the case to which the Sudakov lemma applies, where the stopped process 
$(\tau, X(\tau))$ has the densities ${\bar f}(x, t, \theta)$, and these  can be regarded as the same functions in (\ref{eq:1}) 
(cf.\ e.g.\ Magiera (1974)).

A typical statistical problem is the unbiased estimation of a given $m$-dimensional real vector valued function $h(\theta) = (h_1(\theta), \dots, h_m(\theta))^{\tp}$ of the parameter $\theta$ by using observations on the 
$X(t),\ t\in T$. An estimating procedure for $h(\theta)$ is defined by a pair $(\tau, Z(\tau, X(\tau)))$, where $\tau$ is a stopping variable and $Z(\tau, X(\tau)) = (Z_1(\tau, X(\tau)), \dots, Z_m(\tau, X(\tau)))^{\tp}$ is a ${\mathbb R}^m$-valued function defined on $T \times E$, which is an unbiased estimator of $h(\theta)$, i.e., 
$E_\theta[Z(\tau, X(\tau))] = h(\theta)$.

Let us look at all estimation procedures satisfying the following regularity conditions:
(i)\ The $h(\theta)$ gives a smooth one-to-one transformation from $\Theta$ to $H = h(\Theta)$ in the sense
\begin{align*}
\textrm{rank}\ C_{i \alpha}(\theta) = 
\textrm{rank}\ \frac{\partial h_\alpha}{\partial \theta^i} = 
m,\ \ \forall \theta \in \Theta.
\end{align*}
(ii)\ $E_\theta[\tau] < \infty,\ 
|E_\theta[Z_\alpha(\tau, X(\tau))Z_\beta(\tau, X(\tau))]| < \infty,\ \forall \theta \in \Theta,\ \alpha, \beta = 1, \dots, m$, and the relation $E_\theta[Z(\tau, X(\tau))] = h(\theta)$ 
can be differentiated with respect to $\theta$ under the expectation sign.

Then we have the so-called Cram\'er-Rao inequality for the covariance matrix of the unbiased estimators.

\begin{proposition}
If $Z(\tau, X(\tau))$ is an unbiased estimator of $h(\theta)$, then the the covariance matrix of $Z(\tau, X(\tau))$ is bounded below as 
\begin{align}
\label{eq:3}
E_\theta[(Z(\tau, X(\tau)) - h(\theta))(Z(\tau, X(\tau)) - 
h(\theta))^{\tp}] \ge 
{\bar G}(h(\theta))^{-1},
\end{align}
where
\begin{align*}
&{\bar G}(h(\theta)) = [{\bar g}^{\alpha \beta}(h(\theta))],\ \  
{\bar g}^{\alpha \beta}(h(\theta)) = 
E_\theta[\partial^\alpha {\bar l}_\tau \partial^\beta {\bar l}_\tau] = 
C^{\alpha i}C^{\beta j}{\bar g}_{ij}(\theta),\\
&C^{\alpha i}(h) = \frac{\partial \theta^i}{\partial h_\alpha},\ \ 
{\bar g}_{ij}(\theta) = 
E_\theta[\partial_i {\bar l}_\tau \partial_j {\bar l}_\tau],\ 
{\bar l}_\tau = \log {\bar f}(X(\tau), \tau, h(\theta)),\quad 
\partial^\alpha = \frac{\partial}{\partial h_\alpha},\ \partial_i = 
\frac{\partial}{\partial \theta^i}, 
\end{align*}
and for symmetric matrices $A$ and $B$, the inequality $A \ge B$ implies that $A - B$ is positive semi-definite. The equality in (\ref{eq:3}) holds if and only if $Z(\tau, X(\tau))$ can be represented almost everywhere as
\begin{align}
\label{eq:4}
Z(\tau, X(\tau)) = h(\theta) + {\bar G}(h(\theta))^{-1}
\partial^h {\bar l}(X(\tau), \tau, h(\theta)),\quad 
\partial^h {\bar l} = \bigg(\frac{\partial {\bar l}}{\partial h_1}, \dots, 
\frac{\partial {\bar l}}{\partial h_m} \bigg).
\end{align}
\end{proposition}
\

The condition (\ref{eq:4}) is also written as
\begin{align*}
\partial^h {\bar l}(X(\tau), \tau, h(\theta)) = {\bar G}(h(\theta))
(Z(\tau, X(\tau)) - h(\theta)),
\end{align*}
or in component form
\begin{align*}
\partial^\alpha {\bar l}(X(\tau), \tau, h(\theta)) = 
{\bar g}^{\alpha \beta}(h(\theta))
Z_\beta(\tau, X(\tau)) - {\bar g}^{\alpha \beta}(h(\theta))h_\beta(\theta) = k^\alpha (X(\tau), \tau, h(\theta)).
\end{align*}
The above is a partial differential equation for 
${\bar l}(X(\tau), \tau, h(\theta))$, of which integrability condition $\partial^\beta k^\alpha = 
\partial^\alpha k^\beta$ is
\begin{align*}
[\partial^\beta {\bar g}^{\alpha \gamma}(h(\theta)) - 
\partial^\alpha {\bar g}^{\beta \gamma}(h(\theta))]
(Z_\gamma(\tau, X(\tau)) - h_\gamma(\theta)) = 0,
\end{align*}
and hence the requirement for integrability is
\begin{align*}
\partial^\beta {\bar g}^{\alpha \gamma}(h) =  
\partial^\alpha {\bar g}^{\beta \gamma}(h)\ \Leftrightarrow\ 
\exists {\bar \phi}(h)\ \textrm{smooth and convex in $h$ such that}\ 
{\bar g}^{\alpha \beta}(h) = \partial^\alpha \partial^\beta 
{\bar \phi}(h).
\end{align*}
In this case, the log likelihood function ${\bar l}(x, t, h(\theta))$ is expressed as
\begin{align}
\label{eq:5}
&{\bar l}(x, t, h(\theta)) = {\bar l}(x, t) +
{\bar \xi}^\alpha z_\alpha(t) - {\bar \psi}({\bar \xi}),\\ 
&{\bar \xi}^\alpha = {\bar \xi}^\alpha(h) = \partial^\alpha {\bar \phi}(h),\ \ 
{\bar \psi}({\bar \xi}) = {\bar \xi}^\alpha h_\alpha - {\bar \phi}(h), \nonumber
\end{align}
which implies that $\{{\bar f}(x, t, \theta)\}$ must be a f.r.m.\   exponential family with ${\bar \xi}$ and $h$ the natual and the expectation parameters, respectively.

Suppose that the original $\{f(x, t, \theta)\}$ is not a f.r.m.\  exponential family, then clearly $\{{\bar f}(x, t, \theta)\}$ is not a f.r.m.\ exponential one, either. 
Hence we can restrict attention to the case (\ref{eq:2}), 
when considering the attainment of the lower bound given by Proposition 1.1. 
However we should remark that it is only a necessary condition for the attainment of the lower bound. In fact even in the f.r.m.\ exponential 
family, only some restricted cases can exactly attain the lower bound due to the problem of the ``overshooting'' at the efficient stopping times 
(see Ghosh (1987)).

\section{Conformal transformation of statistical manifold}

Let $M = \{f(x, 1, \theta)\ |\ \theta \in \Theta\}$ be an original family of probability densities of unit time, 
where $\Theta$ is homeomorphic to ${\mathbb R}^m$. 
The family $M$ can be regarded as a statistical manifold, where the 
$m$-dimensiomal vector parameter $\theta$ serves as a coordinate system to specify a point, that is, a density $f(x, 1, \theta)\in M$. 
The geometry of $M$ is determined by the following two tensor quantities 
(cf.\ Amari (1985), Amari and Nagaoka (2000))
\begin{align*}
g_{ij}(\theta) = E_\theta[\partial_il_1\partial_jl_1],\quad
T_{ijk}(\theta) = E_\theta[\partial_il_1\partial_jl_1\partial_kl_1],\ \  
l_1 = \log f(x, 1, \theta),\ \ \partial_i = 
\frac{\partial}{\partial \theta^i}.
\end{align*}
The first is the Fisher information metric and the second is called the skewness tensor. 
One parameter family of affine connections named the $\alpha$-connection is defined by
\begin{align*}
\Gamma_{ijk}^{(\alpha)}(\theta) = 
E_\theta[\partial_i\partial_jl_1\partial_kl_1] + 
\frac{1-\alpha}{2}T_{ijk}(\theta),
\end{align*}
and then the $\alpha$-Riemann-Christoffel curvature tensor is given by
\begin{align*}
R_{ijkl}^{(\alpha)}(\theta) = \partial_i\Gamma_{jkl}^{(\alpha)} - 
\partial_j\Gamma_{ikl}^{(\alpha)} + 
g^{rs}(\Gamma_{ikr}^{(\alpha)}\Gamma_{jsl}^{(\alpha)} - 
\Gamma_{jkr}^{(\alpha)}\Gamma_{isl}^{(\alpha)}).
\end{align*}
The $\alpha$- and the $(-\alpha)$-connections are mutually dual
\begin{align*}
\partial_ig_{jk} = \Gamma_{ijk}^{(\alpha)} + 
\Gamma_{ikj}^{(-\alpha)},
\end{align*}
and the $\pm \alpha$-RC curvature tensors are in the dual relation
\begin{align*}
R_{ijkl}^{(\alpha)} = - R_{ijlk}^{(-\alpha)}.
\end{align*}

Let 
${\bar M} = \{{\bar f}(x, t, \theta)\ |\ \theta \in \Theta\}$ be an $m$-dimensional extended statistical manifold under a sequential statistical procedure.
From the Wald identity the metric and the skewness tensors of 
${\bar M}$ are given by (see Akahira and Takeuchi (1989))
\begin{align}
\label{eq:6}
&{\bar g}_{ij}(\theta) = \nu g_{ij},\quad 
{\bar T}_{ijk}(\theta) = \nu[T_{ijk} + 3g_{(ij}s_{k)}],\\
&\nu(\theta) = E_\theta[\tau],\quad 
s_k(\theta) = \partial_k \log \nu(\theta),\quad 
3g_{(ij}s_{k)} = g_{ij}s_{k} + g_{jk}s_{i} + g_{ki}s_{j}.
\end{align}
These relations show that a sequential statistical procedure induces a conformal transformation $M \mapsto {\bar M}$ by the gauge function $\nu(\theta) > 0$. 
The conformal transformation of a Riemannian manifold implies that the 
manifold is expanded or contracted isotropically but that an expansion rate depends on each point. Our transformation is a statistical counterpart of this one. 
A conformal transformation changes the $\alpha$-connection into
\begin{align}
\label{eq:8}
&{\bar \Gamma}_{ijk}^{(\alpha)} = 
\nu[\Gamma_{ijk}^{(\alpha)} + \frac{1-\alpha}{2}(g_{ki}s_j + g_{kj}s_i)
- \frac{1+\alpha}{2}g_{ij}s_k],\\
&{\bar \Gamma}_{ij}^{(\alpha)k} = {\bar \Gamma}_{ijl}^{(\alpha)}{\bar g}^{lk} = 
\Gamma_{ij}^{(\alpha)k} + \frac{1-\alpha}{2}(\delta_i^ks_j + \delta_j^ks_i)
- \frac{1+\alpha}{2}g_{ij}s_lg^{lk},
\end{align}
This is obtained by substituting (\ref{eq:6}) into
\begin{align*}
{\bar \Gamma}_{ijk}^{(\alpha)} = {\bar \Gamma}_{ijk}^{(0)} - 
\frac{\alpha}{2}{\bar T}_{ijk},
\end{align*}
and by noting that ${\bar \Gamma}_{ijk}^{(0)}$ is the conformal change of 
the Riemannian connection $\Gamma_{ijk}^{(0)}$.
Then a conformal transformation changes the $\alpha$-RC curvature tensor into
\begin{align}
\label{eq:10}
&{\bar R}_{ijkl}^{(\alpha)} = 
\nu[R_{ijkl}^{(\alpha)} - g_{il}s_{jk}^{(\alpha)} + g_{jl}s_{ik}^{(\alpha)}
- g_{jk}s_{il}^{(-\alpha)} + g_{ik}s_{jl}^{(-\alpha)}],\\
&{\bar R}_{ijk}^{(\alpha)l} = {\bar R}_{ijkr}^{(\alpha)}{\bar g}^{rl} = 
R_{ijk}^{(\alpha)l} - \delta_i^ls_{jk}^{(\alpha)} + 
\delta_j^ls_{ik}^{(\alpha)}
- g_{jk}s_{ir}^{(-\alpha)}g^{rl} + g_{ik}s_{jr}^{(-\alpha)}g^{rl},\\
&s_{ij}^{(\alpha)} = \frac{1-\alpha}{2}
[\nabla_i^{(\alpha)}s_j - \frac{1-\alpha}{2}s_is_j + 
\frac{1+\alpha}{4}g_{ij}s_ks_lg^{kl}],\quad 
\nabla_i^{(\alpha)}s_j = \partial_is_j - \Gamma_{ij}^{(\alpha)k}s_k.
\end{align}
This is obtained by substituting (\ref{eq:8}) into
\begin{align*}
{\bar R}_{ijkl}^{(\alpha)} = \partial_i{\bar \Gamma}_{jkl}^{(\alpha)} - 
\partial_j{\bar \Gamma}_{ikl}^{(\alpha)} + 
{\bar g}^{rs}({\bar \Gamma}_{ikr}^{(\alpha)}{\bar \Gamma}_{jsl}^{(\alpha)} - {\bar \Gamma}_{jkr}^{(\alpha)}{\bar \Gamma}_{isl}^{(\alpha)}).
\end{align*}
We note that under a conformal transformation the mutual duality of $\pm\alpha$-connections is preserved 
\begin{align}
\partial_i{\bar g}_{jk} = {\bar \Gamma}_{ijk}^{(\alpha)} + 
{\bar \Gamma}_{ikj}^{(-\alpha)},
\end{align}
and also the dual relation of the $\pm \alpha$-RC curvature tensors is preserved
\begin{align}
\label{eq:14}
{\bar R}_{ijkl}^{(\alpha)} = - {\bar R}_{ijlk}^{(-\alpha)}.
\end{align}
These are confirmed by the direct calculations with (\ref{eq:6}), 
(\ref{eq:8}) and (\ref{eq:10}).

One of the concerns about the conformal transformation is whether a given manifold can be transformed into a desirable space in some sense. 
The main objective from the geometrical viewpoint is the flatness or the straightness, and it has been investigated usually in terms of the Riemannian connection. From the statistical viewpoint, the main objective is the flatness or the straightness in terms of the mutually dual $\pm 1$-connections. 
Thus we say that a statistical manifold $M$ is 
\textit{conformally mixture (exponential) flat} when there exists a gauge function $\nu(\theta) > 0$ such that 
${\bar R}_{ijk}^{(-1)l} = 0\ ({\bar R}_{ijk}^{(1)l} = 0)$ holds. 
Note that by (\ref{eq:14}) $M$ is conformally mixture flat if and only if $M$ is conformally exponential flat.

In view of these observations and also the work of Okamoto (1988), we  introduce the set of $(-1)$-Weyl-Schouten curvature tensors as follows.
\begin{definition} 
\begin{align}
&W_{ijk}^{(-1)l}(\theta) = R_{ijk}^{(-1)l} - 
\frac{1}{m-1}(\delta_i^lR_{jk}^{(-1)} - \delta_j^lR_{ik}^{(-1)}),\\
&W_{ijk}^{(-1)}(\theta) = \frac{1}{m-1}(\nabla_i^{(-1)}R_{jk}^{(-1)} - 
\nabla_j^{(-1)}R_{ik}^{(-1)}),\\
&W_{ij}^{(-1)}(\theta) = R_{ij}^{(-1)} - R_{ji}^{(-1)},\\
&R_{ij}^{(-1)} = R_{lij}^{(-1)l},\quad 
\nabla_i^{(-1)}R_{jk}^{(-1)} = \partial_iR_{jk}^{(-1)} - 
\Gamma_{ij}^{(-1)l}R_{lk}^{(-1)} - \Gamma_{ik}^{(-1)l}R_{jl}^{(-1)}.
\end{align}
\end{definition}
\

From (11) and (12) we have
\begin{align*}
&{\bar R}_{ijk}^{(-1)l} = R_{ijk}^{(-1)l} - \delta_i^ls_{jk}^{(-1)}
+ \delta_j^ls_{ik}^{(-1)},\quad 
s_{jk}^{(-1)} = \nabla_j^{(-1)}s_k - s_js_k\\
&\Rightarrow\ {\bar R}_{jk}^{(-1)} = R_{jk}^{(-1)} - (m-1)s_{jk}^{(-1)}\\
&\Rightarrow\ s_{jk}^{(-1)} = 
-\frac{1}{m-1}({\bar R}_{jk}^{(-1)} - R_{jk}^{(-1)})\\
&\Rightarrow\ {\bar W}_{ijk}^{(-1)l} = W_{ijk}^{(-1)l},\quad
{\bar W}_{ijk}^{(-1)} = W_{ijk}^{(-1)} + W_{ijk}^{(-1)l}s_l,\quad  
{\bar W}_{ij}^{(-1)} = W_{ij}^{(-1)}.
\end{align*}
For the case $m = 2$ we can also directly check  
$W_{ijk}^{(-1)l} \equiv 0$, and hence  
${\bar W}_{ijk}^{(-1)} = W_{ijk}^{(-1)}$.
Then we obtain the following result as to the conditions for the conformal 
mixture (exponential) flatness.

\begin{theorem}
A statistical manifold $M$ is conformally mixture flat (or equivalently exponential flat) if and only if
\begin{quote}
(i)\ $W_{ijk}^{(-1)l} = 0$ when $m = \dim M \ge 3$.\\ 
(ii)\ $W_{ijk}^{(-1)} = 0$ and $W_{ij}^{(-1)} = 0$ when $m = \dim M = 2.$
\end{quote}
For the sake of simplicity we hereafter express the notion such as conformally mixture (or equivalently exponential) flat as conformally $m(e)$-flat.
\end{theorem}
\

\begin{proof}
Consider the relation 
\begin{align*}
s_{jk}^{(-1)} = \nabla_j^{(-1)}s_k - s_js_k = 
-\frac{1}{m-1}({\bar R}_{jk}^{(-1)} - R_{jk}^{(-1)}).
\end{align*}
When ${\bar R}_{jk}^{(-1)} = 0$ we note by the integrability condition that 
\begin{align*}
&\exists s_k\ \textrm{such that}\ \nabla_j^{(-1)}s_k - s_js_k = 
\frac{1}{m-1}R_{jk}^{(-1)}\\
&\Leftrightarrow\ \nabla_i^{(-1)}\nabla_j^{(-1)}s_k - 
\nabla_j^{(-1)}\nabla_i^{(-1)}s_k = -R_{ijk}^{(-1)l}s_l\\
&\Leftrightarrow\ W_{ijk}^{(-1)l}s_l + W_{ijk}^{(-1)} = 0.
\end{align*}

We first prove the necessity. Suppose that $M$ is conformally 
$m$($e$)-flat. 
Then from ${\bar R}_{ijk}^{(-1)l}, {\bar R}_{jk}^{(-1)} = 0$, 
when $m \ge 3$ we have $W_{ijk}^{(-1)l} = {\bar W}_{ijk}^{(-1)l} = 0$. 
When $m = 2$ since $W_{ijk}^{(-1)l} \equiv 0$ we have 
$W_{ijk}^{(-1)} = {\bar W}_{ijk}^{(-1)} = 0$, 
and since there exists a log gauge function $s = \log \nu$ we have $W_{ij}^{(-1)} = {\bar W}_{ij}^{(-1)} = 0$.

We next prove the sufficiency when $m \ge 3$. 
From the Bianchi's second identity (cf.\ Schouten (1954), p.147)
\begin{align*}
&\nabla_l^{(-1)}R_{ijk}^{(-1)l} + \nabla_j^{(-1)}R_{lik}^{(-1)l}
+ \nabla_i^{(-1)}R_{jlk}^{(-1)l} = 0\\
&\Rightarrow\ \nabla_l^{(-1)}R_{ijk}^{(-1)l} = 
\nabla_i^{(-1)}R_{jk}^{(-1)} - \nabla_j^{(-1)}R_{ik}^{(-1)},
\end{align*}
and from $W_{ijk}^{(-1)l} = 0$ we have
\begin{align*}
\nabla_l^{(-1)}W_{ijk}^{(-1)l} = \nabla_l^{(-1)}R_{ijk}^{(-1)l} - 
\frac{1}{m-1}(\nabla_i^{(-1)}R_{jk}^{(-1)} - \nabla_j^{(-1)}R_{ik}^{(-1)})
= (m-2)W_{ijk}^{(-1)} = 0,
\end{align*}
so that $W_{ijk}^{(-1)} = 0$. Then as noted before, there exists a covariant vector field $s_k$ such that 
$\nabla_j^{(-1)}s_k - s_js_k = \frac{1}{m-1}R_{jk}^{(-1)}$. From the Bianchi's first identity (cf.\ Schouten (1954), p.144)
\begin{align*}
R_{ijk}^{(\pm 1)i} + R_{kij}^{(\pm 1)i} + R_{jki}^{(\pm 1)i} = 0\ 
\Rightarrow\ 
R_{jk}^{(\pm 1)} - R_{kj}^{(\pm 1)} + R_{jkil}^{(\pm 1)}g^{li} = 0,
\end{align*}
and from the duality of $R_{ijkl}^{(\pm 1)}$
\begin{align*}
R_{ijkl}^{(1)} + R_{ijkl}^{(-1)} = -R_{ijlk}^{(1)} - R_{ijlk}^{(-1)},
\end{align*}
we have
\begin{align*}
R_{jk}^{(1)} - R_{kj}^{(1)} + R_{jk}^{(-1)} - R_{kj}^{(-1)} = 0.
\end{align*}
On the other hand, from $W_{ijk}^{(-1)l} = 0$ we have
\begin{align*}
&R_{jilk}^{(1)} = R_{ijkl}^{(-1)} = 
\frac{1}{m-1}(g_{il}R_{jk}^{(-1)} - g_{jl}R_{ik}^{(-1)})\\
&\Rightarrow\ R_{il}^{(1)} = R_{jilk}^{(1)}g^{jk} = 
\frac{1}{m-1}(g_{il}R^{(-1)} - R_{il}^{(-1)}),\quad 
R^{(-1)} = R_{jk}^{(-1)}g^{jk}.
\end{align*}
By combining these two relations we obtain
\begin{align*}
\frac{m-2}{m-1}(R_{jk}^{(-1)} - R_{kj}^{(-1)}) = 0\ 
\Rightarrow\ R_{jk}^{(-1)} - R_{kj}^{(-1)} = 0\ 
\Rightarrow\ \partial_js_k - \partial_ks_j = 0,
\end{align*}
and hence there exists a log gauge function $s = \log \nu$ such that 
$s_k = \partial_k s$.

Finally we prove the sufficiency when $m = 2$. From 
$W_{ijk}^{(-1)l} \equiv 0,\ W_{ijk}^{(-1)} = 0$, as noted above,  
there exists a covariant vector field $s_k$ such that 
$\nabla_j^{(-1)}s_k - s_js_k = \frac{1}{m-1}R_{jk}^{(-1)}$. Then from
\begin{align*} 
W_{jk}^{(-1)} = R_{jk}^{(-1)} - R_{kj}^{(-1)} = 0\ 
\Rightarrow\ \partial_js_k - \partial_ks_j = 0,
\end{align*}
there exists a log gauge function $s = \log \nu$ such that 
$s_k = \partial_k s$.

This completes the proof of the theorem.
\end{proof}
\

We further investigate the implications of Theorem 3.1. 
Suppose that the $(-1)$-RC curvature tensor of a statistical manifold 
$M$ is expressed as
\begin{align}
\label{eq:19}
R_{ijkl}^{(-1)} = \lambda(g_{jk}g_{il} - g_{ik}g_{jl}),
\end{align}
where $\lambda$ is constant on $M$. 
In this case we have
\begin{align}
\label{eq:20}
&R_{ijkl}^{(-1)} = -R_{ijlk}^{(-1)} = R_{ijkl}^{(1)}= R_{jilk}^{(-1)}.
\end{align}
A statistical manifold $M$ satisfying 
(\ref{eq:19}) is said to be 
\textit{a space of constant mixture (exponential) curvature}, and $M$ satisfying (\ref{eq:20}) is said to be 
\textit{conjugate mixture (exponential) symmetric}. 
The conjugate $m$($e$)-symmetry is the special notion of the 
conjugate $\pm\alpha$-symmetry introduced by Lauritzen (1987). 
By definition we know
\begin{align*}
\textrm{$M$ is a space of constant $m$($e$)-curvature}\ 
\Rightarrow\ 
\textrm{$M$ is conjugate $m$($e$)-symmetric.}
\end{align*}

The connections among these notions are summarized in the following theorem.

\begin{theorem}
For the conformal $m$($e$)-flatness of a statistical manifold $M$, the following relations hold.
\begin{quote}
(i)\ $M$ is conjugate $m$($e$)-symmetric and is conformally $m$($e$)-flat if and only if $M$ is a space of constant $m$($e$)-curvature.\\
(ii)\ A f.r.m.\ exponential family $M_e$ is always conformally 
$m$($e$)-flat. 
\end{quote}
\end{theorem}
\

\begin{proof}
We first prove the sufficiency of (i). 
Suppose that $M$ is a space of constant $m$($e$)-curvature. 
From 
\begin{align*}
R_{ijkl}^{(-1)} = \lambda(g_{jk}g_{il} - g_{ik}g_{jl}),
\end{align*}
we have
\begin{align*}
R_{jk}^{(-1)} = R_{ijkl}^{(-1)}g^{il} = (m-1)\lambda g_{jk}.
\end{align*}
Then we obtain
\begin{align*}
&W_{ijkl}^{(-1)} = W_{ijk}^{(-1)r}g_{rl} = 
R_{ijkl}^{(-1)} - \frac{1}{m-1}(g_{il}R_{jk}^{(-1)} - g_{jl}R_{ik}^{(-1)}) = 0,\\
&W_{ijk}^{(-1)} = \lambda(\nabla_i^{(-1)}g_{jk} - 
\nabla_j^{(-1)}g_{ik}) = 0\quad 
(\textrm{since}\ \nabla_i^{(\alpha)}g_{jk} = 
\alpha T_{ijk}),\\
&W_{ij}^{(-1)} = (m-1)\lambda(g_{ij} - g_{ji}) = 0.
\end{align*}

We next prove the necessity of (i). 
Suppose that $M$ is conjugate $m$($e$)-symmetric and is 
conformally $m$($e$)-flat. When $m \ge 3$, we have
\begin{align*}
&R_{jilk}^{(-1)} = R_{ijkl}^{(-1)} = \frac{1}{m-1}(g_{il}R_{jk}^{(-1)} - g_{jl}R_{ik}^{(-1)}),\\
&R_{il}^{(-1)} = R_{jilk}^{(-1)}g^{jk} = 
\frac{1}{m-1}(g_{il}R^{(-1)} - R_{il}^{(-1)})\\
&\Rightarrow\ R_{il}^{(-1)} = \frac{R^{(-1)}}{m}g_{il}\\
&\Rightarrow\ R_{ijkl}^{(-1)} = \frac{R^{(-1)}}{m(m-1)}
(g_{jk}g_{il} - g_{ik}g_{jl})\\
&\Rightarrow\ R_{ijk}^{(-1)l} = \rho
(g_{jk}\delta_i^l - g_{ik}\delta_j^l),\quad 
\rho = \frac{R^{(-1)}}{m(m-1)}.
\end{align*}
By substituting this expression into the Bianchi's second identity, we have
\begin{align*}
&\nabla_r^{(-1)}R_{ijk}^{(-1)l} + \nabla_j^{(-1)}R_{rik}^{(-1)l} + 
\nabla_i^{(-1)}R_{jrk}^{(-1)l} = 0\\
&\Rightarrow\ 
\nabla_r^{(-1)}\rho(g_{jk}\delta_i^l - g_{ik}\delta_j^l) + 
\nabla_j^{(-1)}\rho(g_{ik}\delta_r^l - g_{rk}\delta_i^l) + 
\nabla_i^{(-1)}\rho(g_{rk}\delta_j^l - g_{jk}\delta_r^l) = 0\\ 
&\qquad (\textrm{since}\ \nabla_i^{(-1)}g_{jk} = -T_{ijk}\   
\textrm{is symmetric in $i, j, k$})\\
&\Rightarrow\ (m-1)(m-2)\nabla_r^{(-1)}\rho = 
(m-1)(m-2)\partial_r\rho = 0,
\end{align*}
so that $\rho$ is constant on $M$. 
When $m = 2$ we have
\begin{align*}
&\nabla_i^{(-1)}R_{jk}^{(-1)} = \nabla_j^{(-1)}R_{ik}^{(-1)},\quad 
R_{jk}^{(-1)} = \frac{R^{(-1)}}{2}g_{jk}\\
&\Rightarrow\ \partial_iR^{(-1)}g_{jk} = \partial_jR^{(-1)}g_{ik}\quad 
(\textrm{since}\ \nabla_i^{(-1)}g_{jk} = -T_{ijk})\\ 
&\Rightarrow\ \partial_iR^{(-1)} = 0,
\end{align*}
and again $R^{(-1)}$ is constant on $M$. 
This completes the proof of (i).

Since $M_e$ is a space of zero $m$($e$)-curvature $R_{ijkl}^{(\pm 1)} = 0$, (ii) is obtained from (i). 
\end{proof}
\

Figure 1 illustrates the relations among several notions in Theorem 3.2.

\begin{figure}[htbp]
\begin{center}
\includegraphics[width=10cm,height=7cm]{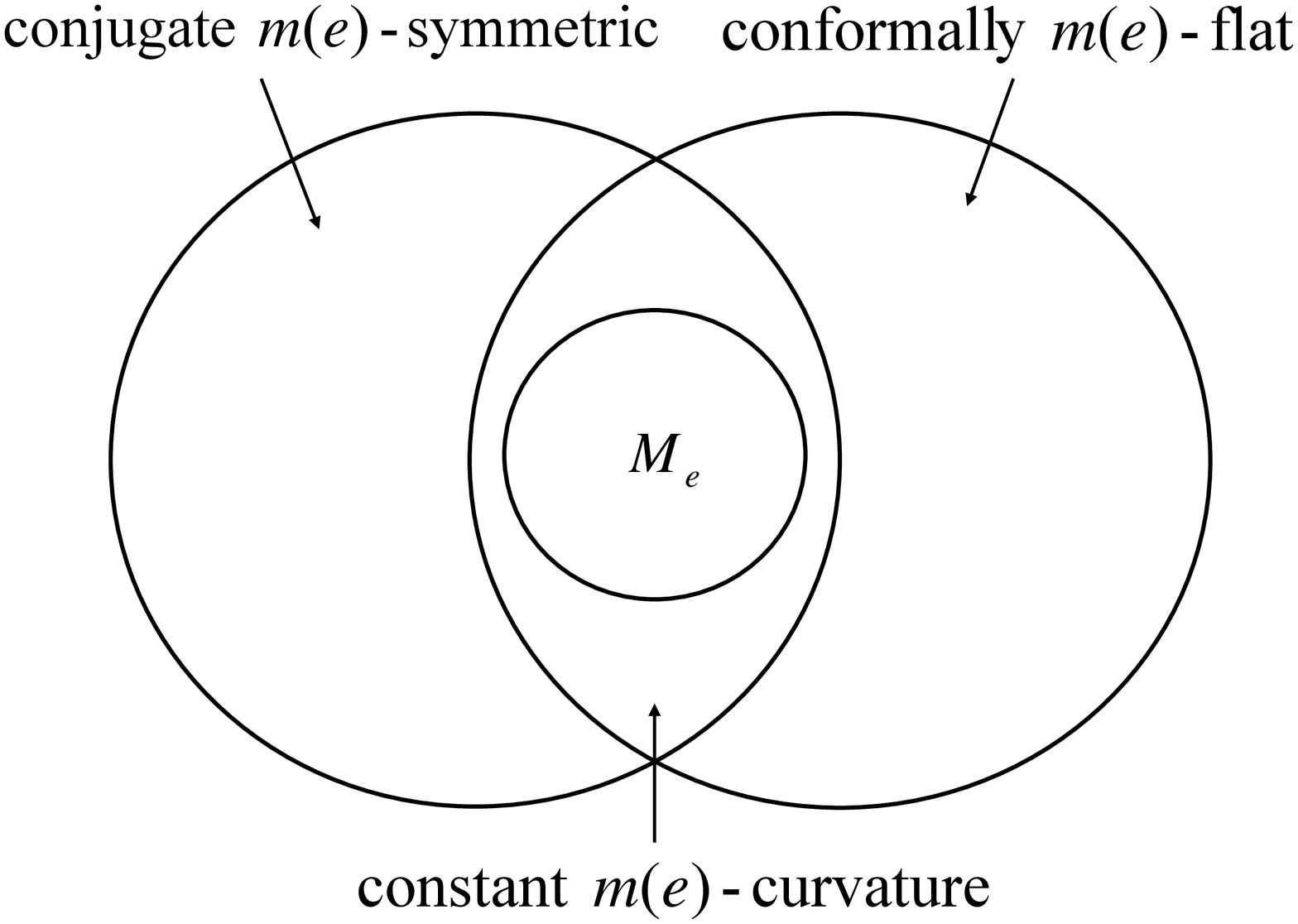}
\vspace*{-7mm}   
\caption{Relations among several notions on $M$}
\label{fig:1-1}
\end{center}
\end{figure}

\section{Conformal geometry of exponential family}

Based on Theorem 3.2 (ii), we seek a concrete conformal transformation 
$M_e \mapsto {\bar M}_e$ such that ${\bar M}_e$ is $\pm 1$-flat.  
When $M_e$ is a f.r.m.\ exponential family, it is $\pm 1$-flat, i.e.,  $R_{ijkl}^{(\pm 1)} = 0$, in which the natural parameter $\theta$ and the expectation parameter $\eta$ provide the $\pm 1$-affine coordinate systems of $M_e$ in the sense (cf.\ Amari (1985), Amari and Nagaoka (2000))
\begin{align*}
&\Gamma_{ijk}^{(1)}(\theta) = 
E_\theta[\partial_i\partial_jl_1\partial_kl_1] = 0,\quad 
\partial_i = \frac{\partial}{\partial \theta^i},\\
&\Gamma^{(-1)ijk}(\eta) = 
E_\eta[(\partial^i\partial^jl_1 + \partial^il_1\partial^jl_1)\partial^kl_1] = 0,\quad \partial^i = \frac{\partial}{\partial \eta_i},
\end{align*}
and there exist two potential functions $\psi(\theta)$ and $\phi(\eta)$ such that 
\begin{align*}
\theta^i = \partial^i \phi(\eta),\ \  
\eta_i = \partial_i \psi(\theta),\quad 
g_{ij}(\theta) = \partial_i\partial_j\psi(\theta),\ \  
g^{ij}(\eta) = \partial^i\partial^j\phi(\eta),\quad 
\psi(\theta) + \phi(\eta) - \theta^i\eta_i = 0.
\end{align*}

By the formula (\ref{eq:8}) a conformal transformation 
$M_e \mapsto {\bar M}_e$ with gauge function $\nu(\eta) > 0$ changes $\Gamma^{(-1)ijk}(\eta)$ into
\begin{align*}
&{\bar \Gamma}^{(-1)ijk}(\eta) = \nu[g^{ki}s^j + g^{kj}s^i],\quad 
g^{ki} = E_\eta[\partial^kl_1\partial^il_1],\ \ s^j = \partial^j\log \nu.
\end{align*}
We consider a coordinate transformation from $\eta$ to $h$ which will provide a $(-1)$-affine coordinate system of ${\bar M}_e$. 
The $(-1)$-connection ${\bar \Gamma}^{(-1)ijk}(\eta)$ transforms to
\begin{align*}
{\bar \Gamma}^{(-1)\alpha \beta \gamma}(h) &= B_i^\alpha B_j^\beta B_k^\gamma 
{\bar \Gamma}^{(-1)ijk}(\eta) + {\bar g}^{jk}B_k^\gamma \partial^\alpha B_j^\beta,\quad B_i^\alpha = \frac{\partial \eta_i}{\partial h_\alpha},\ \ 
\partial^\alpha = \frac{\partial}{\partial h_\alpha},\\
&= \nu B_k^\gamma[B_i^\alpha B_j^\beta (g^{ki}s^j + g^{kj}s^i) + g^{jk}\partial^\alpha B_j^\beta],
\end{align*}
and hence
\begin{align*}
&\textrm{${\bar M}_e$ is $m$-flat}\ \Leftrightarrow\ 
{\bar \Gamma}^{(-1)\alpha \beta \gamma}(h) = 0\\
&\Leftrightarrow\ 
\exists h,\ \nu\ \textrm{such that}\ B_i^\alpha B_j^\beta (g^{ki}s^j + g^{kj}s^i) + g^{jk}\partial^\alpha B_j^\beta = 0\\
&\Leftrightarrow\ \exists h,\ \nu\ \textrm{such that}\ 
\partial^iC_\alpha^j - s^iC_\alpha^j - s^jC_\alpha^i = 0,\ 
s^{(-1)ij} = \partial^is^j - s^is^j = 0,\ C_\alpha^i = 
\frac{\partial h_\alpha}{\partial \eta_i},
\end{align*}
where $s^{(-1)ij} = 0$ is the integrability condition of the first equation on the right-hand side. 
We can solve the above two partial differential equations for 
$s(\eta) = \log \nu(\eta)$ and 
$h_\alpha(\eta)$ as shown in the following theorem.

\begin{theorem}
When a statistical manifold $M_e$ is an $m$-dimensional f.r.m.\  exponential family with $(-1)$-affine coordinate system $\eta$, it is conformally $m$($e$)-flat by the gauge function $\nu(\eta) > 0$ and the new $(-1)$-affine coordinate system $h$ given by
\begin{align}
\label{eq:21}
\nu(\eta) = \frac{1}{|c^0 + c^i\eta_i|},\quad 
h_\alpha = \nu(\eta)(d_\alpha + D_\alpha^i\eta_i),
\end{align}
where $c^0, c^i, d_\alpha, D_\alpha^i\ 
(i, \alpha = 1, \dots, m)$ are constants, 
and \textrm{rank}\ $D_\alpha^i = m$.

The new $1$-affine coordinate system $\xi$, two potential functions 
${\bar \psi}(\xi)$ and ${\bar \phi}(h)$ are respectively given as
\begin{align}
\label{eq:22}
&\xi^\alpha = \partial^\alpha {\bar \phi}(h),\ \ 
h_\alpha = \partial_\alpha {\bar \psi}(\xi),\quad 
{\bar \phi}(h) = \nu(\eta)\phi(\eta),\ \   
{\bar \psi}(\xi) + {\bar \phi}(h) - \xi^\alpha h_\alpha = 0,\\
&{\bar g}_{\alpha \beta}(\xi) = 
\partial_\alpha \partial_\beta {\bar \psi}(\xi),\ \ 
{\bar g}^{\alpha \beta}(h) = 
\partial^\alpha \partial^\beta {\bar \phi}(h).  
\end{align}
\end{theorem}
\

\begin{proof}
We first prove (\ref{eq:21}). By putting $s = -\log r$ we have
\begin{align*}
\partial^is^j - s^is^j = 0\ 
\Leftrightarrow\ \partial^i\partial^j r = 0\ 
\Leftrightarrow\ r = c^0 + c^i\eta_i\ 
\Leftrightarrow\ \nu = \frac{1}{r} = \frac{1}{c^0 + c^i\eta_i},
\end{align*}
and by putting $h_\alpha = e^s y_\alpha$ we have
\begin{align*}
\partial^iC_\alpha^j - s^iC_\alpha^j - s^jC_\alpha^i = 0\ 
\Leftrightarrow\ \partial^i\partial^j y_\alpha = 0\ 
\Leftrightarrow\ y_\alpha = d_\alpha + D_\alpha^i\eta_i\ 
\Leftrightarrow\ h_\alpha = 
\frac{d_\alpha + D_\alpha^i\eta_i}{c^0 + c^i\eta_i}.
\end{align*}

We next prove (22), (23). By the direct calculation we can confirm 
\begin{align*}
{\bar g}^{\alpha \beta}(h) = 
\partial^\alpha \partial^\beta {\bar \phi}(h),\quad 
{\bar \phi}(h) = \nu(\eta)\phi(\eta),
\end{align*}
and then the others are immediately obtained.

This completes the proof of the theorem. 
\end{proof}
\

We remark that $\nu(\eta)$ and $h(\eta)$ in (\ref{eq:21}) cover the general solution and these are the same 
as those given in Winkler and Franz (1979), which were derived from the statistical considerations of the efficient sequential estimators attaining the Cram\'er-Rao bound.

\section{Conformal geometry of curved exponential family}

We first introduce a curved exponential family. A family of probability densities $M_c = \{f_c(x, t, u)\ |\ u\in U\}$ parameterized by an $m$-dimensional vector parameter $u=  (u^1, \dots, u^m)^{\tp}$ is said to be an $(n, m)$-curved exponential family when it is smoothly imbedded in an $n$-dimensional f.r.m.\ exponential family 
$M_e = \{f_e(x, t, \theta)\ |\ \theta\in \Theta\}$ in the sense 
\begin{align}
\label{eq:24}
f_c(x, t, u) = f_e(x, t, \theta(u)) = 
\exp\{\theta^i(u)x_i - \psi(\theta(u))t\},
\end{align}
where $U$ is homeomorphic to 
${\mathbb R}^m (m < n)$ and $\theta(u) = (\theta^1(u), \dots, \theta^n(u))^{\tp}$ is a smooth function of $u$ having a full rank Jacobian matrix. 
We use indices $i, j, k$ and so on to denote quantities in terms of the coordinate system $\theta$ or $\eta$ of $M_e$, and indices $a, b, c$ and so on to denote quantities in terms of the coordinate system $u$ of $M_c$.

For analyzing the 
geometrical properties of $M_c$ imbedded in $M_e$, it is convenient to introduce a new coordinate system $w = (u, v)$ of $M_e$ in the following manner. 
We attach to each point $u\in M_c$ an $(n-m)$-dimensional smooth submanifold $A(u)$ of $M_e$ which transverses $M_c$ at $\theta(u)$ 
or equivalently at $\eta(u)$. We assume that the family $A = \{A(u)\ |\ u\in M_c\}$ fills up at least a neighborhood of $M_c$ in $M_e$, that is, $A$ is a foliation of the tubular neighborhood of $M_c$ in $M_e$. Such an $A(u)$ is called an ancillary submanifold rigging $u$, and $A$ is called an ancillary family rigging $M_c$.

We introduce a coordinate system $v = (v^{m+1}, \dots, v^n)$ to each $A(u)$ such that the origin $v = 0$ is at the intersection of $A(u)$ and $M_c$. Then the combined system
\begin{align*}
w = (w^\alpha) = (u^a, v^\kappa),\ \alpha = 1, \dots, n,\ a = 1, \dots, m,\ \kappa = m+1, \dots, n
\end{align*}
gives a new local coordinate system of $M_e$. 
We use indices $\alpha, \beta, \gamma$ and so on for quantities related to the coordinate system $w$, and indices $\kappa, \lambda, \mu$ and so on for quantities related to the coordinate system $v$.

The basic tensors of $M_e$ are written as
\begin{align*}
g_{\alpha \beta}(w) = g_{ij}(\theta)B_\alpha^iB_\beta^j,\quad 
T_{\alpha \beta \gamma}(w) = T_{ijk}(\theta)B_\alpha^iB_\beta^jB_\gamma^k,\quad 
B_\alpha^i = \frac{\partial \theta^i}{\partial w^\alpha},
\end{align*}
and the $\alpha$-connection is given by 
\begin{align*}
\Gamma_{\beta \gamma \delta}^{(\alpha)}(w) &= 
\Gamma_{ijk}^{(\alpha)}(\theta)B_\beta^iB_\gamma^jB_\delta^k + 
g_{ij}(\theta)B_\delta^i\partial_\beta B_\gamma^j = 
\frac{1-\alpha}{2}T_{\beta \gamma \delta} + 
(\partial_\beta B_\gamma^j)B_{\delta j}\\
&= \Gamma^{(\alpha)ijk}(\eta)B_{\beta i}B_{\gamma j}B_{\delta k} + 
g^{ij}(\eta)B_{\delta i}\partial_\beta B_{\gamma j} = 
-\frac{1+\alpha}{2}T_{\beta \gamma \delta} + 
(\partial_\beta B_{\gamma j})B_\delta^j,\\
B_{\beta i} &= \frac{\partial \eta_i}{\partial w^\beta} = g_{ij}B_\beta^j,
\end{align*}
in the $w$-coordinate system. When we evaluate a quantity $q(u, v)$ on $M_c$, i.e., at $v = 0$, we often denote it by $q(u)$ instead of by $q(u, 0)$ for brevity's sake. 
The metric tensors of $M_c$ and $A(u)$ are given by
\begin{align*}
g_{ab}(u) = B_a^iB_b^jg_{ij} = B_{ai}B_{bj}g^{ij},\quad 
g_{\kappa\lambda}(u) = B_\kappa^iB_\lambda^jg_{ij} = 
B_{\kappa i}B_{\lambda j}g^{ij},
\end{align*}
and then indices can be lowered or uppered by using these metric tensors or their inverses $g^{ab}(u), g^{\kappa \lambda}(u)$.
The $\pm1$-connections of $M_c$ are given by
\begin{align*}
\Gamma_{abc}^{(1)}(u) = (\partial_a B_b^j)B_{cj},\quad 
\Gamma_{abc}^{(-1)}(u) = (\partial_a B_{b j})B_c^j.
\end{align*}
We call $A = \{A(u)\ |\ u\in U\}$ an orthogonal ancillary family when 
$g_{a\kappa}(u) = 0,\ \forall u\in U$, and we assume this property in the following.  
The mixed parts $\Gamma_{ab\kappa}^{(\pm1)}(u)$ play central roles in the evaluation of statistical inferences, which are defined as follows.
\begin{definition}
\begin{align*}
H_{ab\kappa}^{(1)}(u) = \Gamma_{ab\kappa}^{(1)}(u) = 
(\partial_a B_b^j)B_{\kappa j},\quad
H_{ab\kappa}^{(-1)}(u) = \Gamma_{ab\kappa}^{(-1)}(u) = 
(\partial_a B_{bj})B_\kappa^j,
\end{align*}
and we call $H_{ab\kappa}^{(\pm1)}$ the $\pm1$-Euler-Schouten curvature tensors of $M_c$. 
\end{definition}
\

The $\pm1$-RC curvature tensors and the $\pm$-ES curvature tensors of 
$M_c$ are connected by the equations of Gauss 
(cf.\ Schouten (1954), p.266)
\begin{align}
\label{eq:25}
R_{abcd}^{(\pm1)}(u) &= R_{ijkl}^{(\pm1)}B_a^iB_b^jB_c^kB_d^l + 
(H_{ad\kappa}^{(\mp1)}H_{bc\lambda}^{(\pm1)} - 
H_{bd\kappa}^{(\mp1)}H_{ac\lambda}^{(\pm1)})g^{\kappa \lambda}\nonumber \\
&= (H_{ad\kappa}^{(\mp1)}H_{bc\lambda}^{(\pm1)} - 
H_{bd\kappa}^{(\mp1)}H_{ac\lambda}^{(\pm1)})g^{\kappa \lambda}.
\end{align}
Suppose that the $\pm1$-ES curvature tensors of $M_c$ are related as
\begin{align}
\label{eq:26}
H_{ab\kappa}^{(-1)}(u) = \epsilon H_{ab\kappa}^{(1)}(u),
\end{align}
where $\epsilon\ (\neq 0)$ is a constant. 
In this case we have
\begin{align*}
R_{abcd}^{(-1)}(u) = R_{abcd}^{(1)}(u) = 
\epsilon (H_{ad\kappa}^{(1)}H_{bc\lambda}^{(1)} - 
H_{bd\kappa}^{(1)}H_{ac\lambda}^{(1)})g^{\kappa \lambda},
\end{align*}
so that $M_c$ is conjugate $m$($e$)-symmetric. Thus we say that $M_c$ satisfying (\ref{eq:26}) is 
\textit{ES conjugate $m$($e$)-symmetric}. 
Suppose further that the $1$-ES curvature tensor of $M_c$ is written as
\begin{align}
\label{eq:27}
H_{ab\kappa}^{(1)}(u) = H_{\kappa}^{(1)}g_{ab}(u),\quad 
H_{\kappa}^{(1)}(u) = \frac{1}{m}H_{ab\kappa}^{(1)}g^{ab},\quad 
\forall u\in M_c,
\end{align}
where $H_{\kappa}^{(1)}$ is called the mean $1$-ES curvature of 
$M_c$, and $M_c$ satisfying (\ref{eq:27}) 
is said to be 
\textit{totally exponential umbilic (e-umbilic)}.

The implications of these notions are summarized in the following manner.

\begin{theorem}
For an $(n, m)$-curved exponential family $M_c$, the following relation
 holds.
\begin{quote}
Let $m \ge 3$ or $n = m+1$, and suppose that 
$M_c$ is ES conjugate $m$($e$)-symmetric and totally $e$-umbilic. 
Then  
$M_c$ is a space of constant $m$($e$)-curvature.
\end{quote}
\end{theorem}
\

\begin{proof}
Suppose that $M_c$ is ES conjugate $m$($e$)-symmetric and totally $e$-umbilic. Then from 
(\ref{eq:26}) and (\ref{eq:27}) we have
\begin{align*}
R_{abcd}^{(-1)}(u) = R_{abcd}^{(1)}(u) = 
\epsilon H^{(1)2}(g_{ad}g_{bc} - g_{ac}g_{bd}),\quad 
H^{(1)2}(u) = H_\kappa^{(1)}H_\lambda^{(1)}g^{\kappa \lambda}.
\end{align*}
When $m \ge 3$, as noted in the proof of Theorem 3.2, $\epsilon H^{(1)2}$ is constant on $M_c$. When $n = m + 1$, from the equation of Codazzi 
(cf.\ Schouten (1954), p.266)
\begin{align*}
&0 = R_{ijk}^{(1)l}B_a^iB_b^jB_c^kB_l^\kappa = 
\nabla_a^{(1)}H_{bc}^{(1)\kappa} - \nabla_b^{(1)}H_{ac}^{(1)\kappa},\quad 
B_l^\kappa = g_{li}g^{\kappa \lambda}B_{\lambda}^i\\
&\Rightarrow\ 
\epsilon(\nabla_a^{(1)}H^{(1)\kappa})g_{bc} -  
\epsilon(\nabla_b^{(1)}H^{(1)\kappa})g_{ac} = 0\quad 
(\textrm{since}\ \nabla_a^{(1)}g_{bc} = T_{abc})\\
&\Rightarrow\ \epsilon(m - 1)\nabla_a^{(1)}H^{(1)\kappa} = 
\epsilon(m - 1)\partial_aH^{(1)\kappa} = 0\quad 
(\textrm{since}\ 
\textrm{$H^{(1)\kappa}$ is a scalar on $M_c$,\ $\kappa = m+1$}),
\end{align*}
and without loss of generality we can set 
$g_{\kappa \kappa}(u) = 1$, 
so that $\epsilon H^{(1)2} = 
\epsilon H^{(1)\kappa}H^{(1)\kappa}g_{\kappa \kappa}$ 
is again constant on $M_c$. 
\end{proof}
\

We further deal with the case of $n = m + 1$. Suppose that 
$M_c$ satisfies the equations
\begin{align}
\label{eq:28}
B_\kappa^i(u) = k_0(\theta^i(u) - \theta_0^i),\quad
B_{\kappa i}(u) = l_0(\eta_i(u) - \eta_i^0),\quad g_{\kappa \kappa}(u) = 1,
 \ \kappa = m + 1,
\end{align}
where $k_0, l_0$ are non-zero constants and $\theta^i_0, \eta_i^0$ 
are constant vectors. In this case $M_c$ is expressed as 
\begin{align*}
(\theta^i(u) - \theta^i_0)(\eta_i(u) - \eta_i^0) = \frac{1}{k_0l_0},
\end{align*}
and we call $M_c$ satisfying (\ref{eq:28}) 
\textit{a dual quadric hypersurface}. 
In Section 7 it will be shown that the von Mises-Fisher model and the hyperboloid model are the examples of the dual quadric hypersurface. 
The meaning of this hypersurface is described in the following theorem.

\begin{theorem}
For an $(m+1, m)$-curved exponential family $M_c$, the following two conditions are equivalent.
\begin{quote}
(i)\ $M_c$ is a dual quadric hypersurface.\\
(ii)\ $M_c$ is ES congugate $m$($e$)-symmetric and totally $e$-umbilic with constant \\
\quad \ $m$($e$)-curvature $k_0l_0$, and $T_{a\kappa \kappa}(u) = 0$ on $M_c$.
\end{quote}
\end{theorem}
\

\begin{proof}
We first prove $(i) \Rightarrow (ii)$.   
By the definition (\ref{eq:28}) we have
\begin{align*}
&\partial_aB_\kappa^i(u) = 
\Gamma_{a\kappa}^{(1)b}B_{bi} + 
\Gamma_{a\kappa}^{(1)\kappa}B_{\kappa i} = k_0B_a^i(u),\quad 
\partial_aB_{\kappa i}(u) = 
\Gamma_{a\kappa}^{(-1)b}B_{bi} + 
\Gamma_{a\kappa}^{(-1)\kappa}B_{\kappa i} = l_0B_{ai}(u)\\
&\Rightarrow\ 
\Gamma_{a\kappa b}^{(1)}(u) = k_0g_{ab}(u),\quad 
\Gamma_{a\kappa \kappa}^{(1)}(u) = 0,\quad 
\Gamma_{a\kappa b}^{(-1)}(u) = l_0g_{ab}(u),\quad 
\Gamma_{a\kappa \kappa}^{(-1)}(u) = 0.
\end{align*}
On the other hand
\begin{align*}
0 = \partial_bg_{a\kappa}(u) = 
H_{ab\kappa}^{(-1)} + \Gamma_{b\kappa a}^{(1)} = 
H_{ab\kappa}^{(1)} + \Gamma_{b\kappa a}^{(-1)},
\end{align*}
and hence
\begin{align*}
&H_{ab\kappa}^{(-1)}(u) = -k_0g_{ab}(u) = 
\frac{k_0}{l_0}H_{ab\kappa}^{(1)}(u),\quad
H_{ab\kappa}^{(1)}(u) = -l_0g_{ab}(u),\quad 
T_{a\kappa \kappa}(u) = 
\Gamma_{a\kappa \kappa}^{(-1)} - \Gamma_{a\kappa \kappa}^{(1)}= 0\\
&\Rightarrow\ R_{abcd}^{(\pm1)}(u) = 
k_0l_0(g_{ad}g_{bc} - g_{ac}g_{bd}).
\end{align*}

We next prove $(ii) \Rightarrow (i)$. 
By the definitions (\ref{eq:26}) and 
(\ref{eq:27}) we have
\begin{align*}
&R_{abcd}^{(\pm1)}(u) = 
\epsilon H^{(1)2}(g_{ad}g_{bc} - g_{ac}g_{bd})\\
&\Rightarrow \epsilon H^{(1)2} = 
\epsilon H_\kappa^{(1)}H_\kappa^{(1)}g^{\kappa \kappa} = k_0l_0,\quad 
g^{\kappa \kappa}(u) = 1\\
&\Rightarrow\ H_{ab\kappa}^{(1)}(u) = H_\kappa^{(1)}g_{ab}(u) = 
-\sqrt{|k_0l_0/\epsilon|}g_{ab}(u),\quad 
H_{ab\kappa}^{(-1)}(u) = -\sqrt{|\epsilon k_0l_0|}g_{ab}(u).
\end{align*}
On the ther hand 
\begin{align*}
&\Gamma_{b\kappa a}^{(1)}(u) = -H_{ab\kappa}^{(-1)}(u),\quad
\Gamma_{b\kappa a}^{(-1)}(u) = -H_{ab\kappa}^{(1)}(u),\\
&0 = \partial_ag_{\kappa \kappa}(u) = 
\Gamma_{a\kappa \kappa}^{(1)}(u) + \Gamma_{a\kappa \kappa}^{(-1)}(u),\quad 
T_{a\kappa \kappa}(u) = 
\Gamma_{a\kappa \kappa}^{(-1)} - \Gamma_{a\kappa \kappa}^{(1)} = 0,
\end{align*}
and hence
\begin{align*}
&\Gamma_{a\kappa \kappa}^{(1)}(u)= 0,\quad 
\Gamma_{a\kappa \kappa}^{(-1)}(u)= 0\\
&\Rightarrow\ 
\partial_aB_\kappa^i(u) = \Gamma_{a\kappa}^{(1)b}B_b^i(u) + 
\Gamma_{a\kappa}^{(1)\kappa}B_\kappa^i(u) = 
\sqrt{|\epsilon k_0l_0|}B_a^i(u)\\
&\qquad \partial_aB_{\kappa i}(u) = \Gamma_{a\kappa}^{(-1)b}B_{bi}(u) + 
\Gamma_{a\kappa}^{(-1)\kappa}B_{\kappa i}(u) = 
\sqrt{|k_0l_0/\epsilon|}B_{ai}(u)\\
&\Rightarrow\ \partial_a(B_\kappa^i(u) - 
\sqrt{|\epsilon k_0l_0|}\theta^i(u)) = 0,\quad 
\partial_a(B_{\kappa i}(u) - 
\sqrt{|k_0l_0/\epsilon|}\eta_i(u)) = 0\\
&\Rightarrow\ B_\kappa^i(u) = 
\sqrt{|\epsilon k_0l_0|}(\theta^i(u) - \theta^i_0),\quad 
B_{\kappa i}(u) = \sqrt{|k_0l_0/\epsilon|}(\eta_i(u) - \eta_i^0).
\end{align*}

This completes the proof of the theorem.
\end{proof}
\

Figure 2 illustrates the relations among several notions in Theorems 5.1 
and 5.2.

\begin{figure}[htbp]
\begin{center}
\includegraphics[width=10cm,height=7cm]{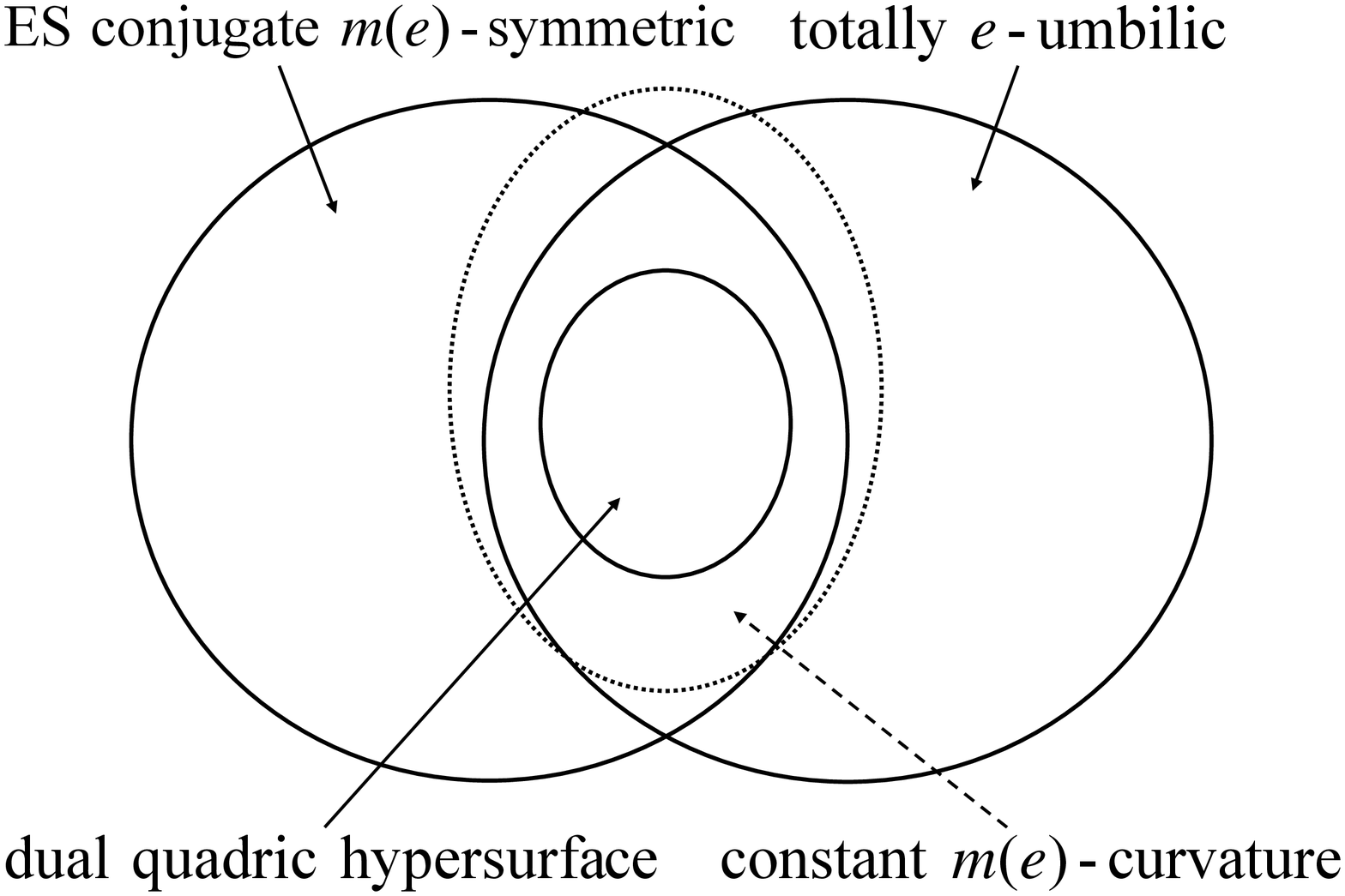}
\vspace*{-7mm}   
\caption{Relations among several notions on $M_c$}
\label{fig:1-1}
\end{center}
\end{figure}

From Theorems 3.2 and 5.2, the dual quadric hypersurface is conformally 
$m(e)$-flat, and we obtain the following result as to its dual structure.

\begin{theorem}
When $M_c$ is an $m$-dimensional dual quadric hypersurface, it is conformally $m(e)$-flat by the log gauge function $s(u) = \log \nu(u)$ and the $(-1)$-affine coordinate system 
${\bar u}^{\bar a},\ {\bar a} = 1, \dots, m$,  
satisfying 
\begin{align}
\label{eq:29}
\partial_a s_b - \Gamma_{ab}^{(-1)c}s_c - s_as_b = k_0l_0g_{ab},\quad
{\bar u}^{\bar a} = \nu(u)D^{{\bar a}i}(\eta_i(u) - \eta_i^0),
\end{align}
where 
$\eta_i^0, D^{{\bar a}i}\ ({\bar a} = 1, \dots, m,\ i = 1, \dots, m+1)$ are constants and  
\textrm{rank}\ $D^{{\bar a}i} = m$.

The $1$-affine coordinate system ${\bar \upsilon}_{\bar a}$, 
two potential functions ${\bar \psi}({\bar \upsilon})$ and 
${\bar \phi}({\bar u})$ of ${\bar M}_c$ are respectively given as 
\begin{align}
\label{eq:30}
&{\bar \upsilon}_{\bar a} = \partial_{\bar a} {\bar \phi}({\bar u}),\ \ 
{\bar u}^{\bar a} = \partial^{\bar a} {\bar \psi}({\bar \upsilon}),\quad
{\bar \phi}({\bar u}) = \frac{\nu(u)}{k_0l_0},\ \ 
{\bar \psi}({\bar \upsilon}) + {\bar \phi}({\bar u}) - 
{\bar \upsilon}_{\bar a}{\bar u}^{\bar a} = 0,\\
&g^{{\bar a}{\bar b}}({\bar \upsilon}) = 
\partial^{\bar a} \partial^{\bar b} {\bar \psi}({\bar \upsilon}),\ \ 
g_{{\bar a}{\bar b}}({\bar u}) = 
\partial_{\bar a} \partial_{\bar b} {\bar \phi}({\bar u}).
\end{align}
\end{theorem}
\

\begin{proof}
We first prove (\ref{eq:29}). 
As noted in the proof of Theorem 3.1, the partial differential equation for 
$s(u) = \log \nu(u)$ is
\begin{align*}
\partial_a s_b - \Gamma_{ab}^{(-1)c}s_c - s_as_b = 
\frac{1}{m-1}R_{ab}^{(\pm1)}.
\end{align*}
When $M_c$ is an $m$-dimensional dual quadric hypersurface, from 
\begin{align*}
R_{abcd}^{(\pm1)} = k_0l_0(g_{ad}g_{bc} - g_{ac}g_{bd})\ \Rightarrow\ 
R_{bc}^{(\pm1)} = R_{abcd}^{(\pm1)}g^{ad} = k_0l_0(m-1)g_{bc},
\end{align*}
we have the first relation of (\ref{eq:29}).

The partial differential equation for ${\bar u} = {\bar u}(u)$ is given as
\begin{align*}
&{\bar \Gamma}_{{\bar a}{\bar b}{\bar c}}^{(-1)}({\bar u}) = 
B_{\bar a}^aB_{\bar b}^bB_{\bar c}^c {\bar \Gamma}_{abc}^{(-1)}(u) + 
{\bar g}_{bc}B_{\bar c}^c \partial_{\bar a}B_{\bar b}^b = 0\\
&\Leftrightarrow\ 
B_{\bar c}^c[B_{\bar a}^aB_{\bar b}^b{\bar \Gamma}_{abc}^{(-1)} + 
g_{bc}\partial_{\bar a}B_{\bar b}^b] = 0,\quad 
{\bar \Gamma}_{abc}^{(-1)} = 
\nu[\Gamma_{abc}^{(-1)} + g_{ac}s_b + g_{bc}s_a]\\
&\Leftrightarrow\ 
\partial_bC_a^{\bar a} - {\bar \Gamma}_{ba}^{(-1)c}C_c^{\bar a} = 0,\quad 
C_a^{\bar a} = \frac{\partial {\bar u}^{\bar a}}{\partial u^a},\quad 
{\bar \Gamma}_{ba}^{(-1)c} = \Gamma_{ba}^{(-1)c} + 
\delta_a^cs_b + \delta_b^cs_a.
\end{align*}
When $M_c$ is a dual quadric hypersurface, 
we can directly show that the above is satisfied by the second relation of (\ref{eq:29}).

We next prove (\ref{eq:30}), (31). By the direct calculation we can confirm
\begin{align*}
g_{{\bar a}{\bar b}}({\bar u}) = 
\partial_{\bar a} \partial_{\bar b} {\bar \phi}({\bar u}),\quad 
{\bar \phi}({\bar u}) = \frac{\nu(u)}{k_0l_0},
\end{align*}
and then the others are immediately obtained.

This completes the proof of the theorem.
\end{proof}
\

We next consider a conformal transformation $M_e \mapsto {\bar M}_e$ by the gauge function $\nu(w) > 0$. 
As shown by (\ref{eq:8}), the $\alpha$-connection in terms of the $w$-coordinate system is changed into
\begin{align*}
&{\bar \Gamma}_{\beta \gamma \delta}^{(\alpha)} = 
\nu[\Gamma_{\beta \gamma \delta}^{(\alpha)} + \frac{1-\alpha}{2}
(g_{\delta \beta}s_\gamma + g_{\delta \gamma}s_\beta)
- \frac{1+\alpha}{2}g_{\beta \gamma}s_\delta].
\end{align*}
Then we can express the change of quantities related to ${\bar M}_c$, that is, the $(-1)$-connection of $M_c$, the $1$-ES curvature of 
$M_c$ and the $(-1)$-ES curvature of $A(u)$ are respectively changed into
\begin{align}
\label{eq:32}
{\bar \Gamma}_{abc}^{(-1)} = 
\nu[\Gamma_{abc}^{(-1)} + g_{ca}s_b + g_{cb}s_a],\quad
{\bar H}_{ab\kappa}^{(1)} = 
\nu[H_{ab\kappa}^{(1)} - g_{ab}s_\kappa],\quad
{\bar H}_{\kappa \lambda a}^{(-1)} = 
\nu H_{\kappa \lambda a}^{(-1)}.
\end{align}
It is also seen that
\begin{align}
\label{eq:33}
{\bar K}_{ab\kappa}^{(1)} = \nu K_{ab\kappa}^{(1)},\quad 
K_{ab\kappa}^{(1)}(u) = H_{ab\kappa}^{(1)} - g_{ab}H_{\kappa}^{(1)},
\end{align}
and we call $K_{ab\kappa}^{(1)}$ 
\textit{the conformal $1$-ES curvature tensor}.

Note that the change of the $(-1)$-connection
\begin{align*}
{\bar \Gamma}_{\alpha \beta}^{(-1)\gamma} = 
\Gamma_{\alpha \beta}^{(-1)\gamma} + 
\delta_\alpha^\gamma s_\beta + \delta_\beta^\gamma s_\alpha
\end{align*}
induces the projective transformation at the same time, which implies that the mixture geodesic is preserved under the transformation (cf.\ Schouten (1954), p.287). 
The effect of constant $m$($e$)-curvature is given in the following theorem.

\begin{theorem}
Suppose that a curved exponential family $M_c$ is a space of constant 
$m$($e$)-curvature. 
Then there exists a conformal transformation $M_e \mapsto {\bar M}_e$ and a coordinate system ${\bar u} = ({\bar u}^{{\bar a}}),\ {\bar a} = 1, \dots, m$, of ${\bar M}_c$ such that the followings hold.
\begin{quote}
(i)\ ${\bar \Gamma}_{{\bar a}{\bar b}{\bar c}}^{(-1)}({\bar u}) = 0,\ \  
\forall {\bar u}\in {\bar M}_c.$\\
(ii)\ If $M_c$ is totally $e$-umbilic, then 
${\bar H}_{{\bar a}{\bar b}\kappa}^{(1)}({\bar u}) = 0,\ \  
\forall {\bar u}\in {\bar M}_c.$
\end{quote}
\end{theorem}
\

\begin{proof}
We first prove (i). When $M_c$ is a space of constant $m$($e$)-curvature, from Theorem 3.2, $M_c$ is conformally $m$($e$)-flat, so that we have
\begin{align*}
{\bar R}_{abcd}^{(-1)}(u) = 0\ \Leftrightarrow\ 
\exists \nu(u) > 0,\ {\bar u} = ({\bar u}^{{\bar a}})\ \textrm{such that}\ 
{\bar \Gamma}_{{\bar a}{\bar b}{\bar c}}^{(-1)}({\bar u}) = 0,\ 
\forall {\bar u}\in {\bar M}_c.
\end{align*}

We next prove (ii). Let us take $s_\kappa(u) = H_\kappa^{(1)}(u)$ 
(see Okamoto, Amari and Takeuchi (1991)).   
Then for the totally $e$-umbilic $M_c$, from (\ref{eq:27}), (\ref{eq:32})  and (\ref{eq:33}) we have
\begin{align*}
{\bar H}_{ab\kappa}^{(1)}(u) = {\bar K}_{ab\kappa}^{(1)}(u) =
\nu K_{ab\kappa}^{(1)}(u) = 0,\ \forall u\in {\bar M}_c
\Rightarrow\ 
{\bar H}_{{\bar a}{\bar b}\kappa}^{(1)}({\bar u}) = 0,\ 
\forall {\bar u}\in {\bar M}_c.
\end{align*}
\end{proof}

\section{Sequential estimation in curved exponential family}

We consider sequential estimations in an $(n, m)$-curved exponential family $M_c$. 
Let $K > 0$ be a parameter cotrolling the average sample size, and let $\nu(\eta) > 0\ (\nu(w) > 0)$ be a smooth gauge function defined on $M_e$ in the $\eta$-($w$-)coordinate system.

We denote by ${\bar X}_t = X(t)/t$ the sample mean up to time $t$. 
It has the same value in the $\eta$-coordinate system, and its value 
in the $w$-coordinate system is denoted by 
${\hat w}_t = ({\hat u}_t, {\hat v}_t) = \eta^{-1}({\bar X}_t)$. 
The random stopping time $\tau$ is assumed to satisfy 
(see Okamoto, Amari and Takeuchi (1991)) 
\begin{align*}
&\tau = K\nu({\hat w}_\tau) + c({\hat u}_\tau) + \varepsilon,\quad 
c(u) = -\frac{1}{2}(\partial_\alpha s_\beta - 
\Gamma_{\alpha \beta}^{(-1)\gamma}s_\gamma - 
s_\alpha s_\beta)g^{\alpha \beta},\\ 
&\varepsilon = O_p(1),\ E_u[\varepsilon] = o(1),\quad 
E_u[\tau] = K\nu(u),\ V_u[\tau] = O(K).
\end{align*}
The term $c$ is due to the bias of ${\hat w}_\tau$ from the true 
$w = (u, 0)$, which is obtained by the requirement 
$E_u[\tau] = K\nu(u)$. The term $\varepsilon$ includes a rounding error 
and the ``overshooting'' at the stopping time $\tau$.

We cite the established results concerning the asymptotics of sequential estimators of $u$ from Okamoto, Amari and Takeuchi (1991).

\begin{proposition}
For a consistent sequential estimator ${\hat u}$ of $u$, the following relations hold.

\noindent
(i)\ The estimator ${\hat u}$ is first-order efficient, that is, 
$\sqrt{K\nu}({\hat u} - u) \to N(0, g^{ab}(u))$ as $K \to \infty$, 
if and only if $A = \{A(u)\}$ is an orthogonal ancillary family.\\

\noindent
(ii)\ The bias-corrected estimator ${\hat u}^*$ of ${\hat u}$ is given by
\begin{align}
\label{eq:34}
{\hat u}^{*a} = {\hat u}^a + 
\frac{1}{2K\nu}\ ^{'}\!\Gamma_{\alpha \beta}^{(-1)a}
g^{\alpha \beta}({\hat u}),\quad
\ ^{'}\!\Gamma_{\alpha \beta}^{(-1)a} = \Gamma_{\alpha \beta}^{(-1)a} + 
\delta_\alpha^a s_\beta + \delta_\beta^a s_\alpha.
\end{align}

\noindent
(iii)\ 
The asymptotic covariance of ${\hat u}^*$ is given by
\begin{align}
\label{eq:35}
E[K\nu({\hat u}^{*a} - u^a)({\hat u}^{*b} - u^b)] = 
g^{ab} + \frac{1}{K\nu}\bigg\{\frac{1}{2} 
(^{'}\!\Gamma_{M_c}^{(-1)})^{2ab} + 
 (^{'}\!H_{M_c}^{(1)})^{2ab} + \frac{1}{2}(H_A^{(-1)})^{2ab} \bigg\} + O(K^{-2}),
\end{align}
where
\begin{align*}
&(^{'}\!\Gamma_{M_c}^{(-1)})^{2ab} = 
(^{'}\!\Gamma_{M_c}^{(-1)})_{cd}^2g^{ac}g^{bd},\quad 
(^{'}\!\Gamma_{M_c}^{(-1)})_{ab}^2 =
\ ^{'}\!\Gamma_{cda}^{(-1)}\ ^{'}\!\Gamma_{efb}^{(-1)}g^{ce}g^{df},\quad 
\ ^{'}\!\Gamma_{abc}^{(-1)} = \Gamma_{abc}^{(-1)} + g_{ca}s_b + g_{cb}s_a,\\
&(^{'}\!H_{M_c}^{(1)})^{2ab} = 
(^{'}\!H_{M_c}^{(1)})_{cd}^2g^{ac}g^{bd},\quad 
(^{'}\!H_{M_c}^{(1)})_{ab}^2 =
\ ^{'}\!H_{ac\kappa}^{(1)}\ ^{'}\!H_{bd\lambda}^{(1)}g^{cd}
g^{\kappa \lambda},\quad \ ^{'}\!H_{ab\kappa}^{(1)} = H_{ab\kappa}^{(1)} - g_{ab}s_\kappa, \\
&(H_A^{(-1)})^{2ab} = 
(H_A^{(-1)})_{cd}^2g^{ac}g^{bd},\quad 
(H_A^{(-1)})_{ab}^2 =
H_{\kappa \lambda a}^{(-1)}H_{\mu \nu b}^{(-1)}
g^{\kappa \mu}g^{\lambda \nu}.
\end{align*}
\end{proposition}
\

Based on Theorem 5.4 we obtain the following result for the possibility of covariance minimization.

\begin{theorem}
Suppose that a curved exponential family $M_c$ is a space of constant 
$m$($e$)-curvature and is totally $e$-umblic. 
Then there exists a conformal transformation $M_e \mapsto {\bar M}_e$ and a coordinate system ${\bar u} = ({\bar u}^{{\bar a}}),\ {\bar a} = 1, \dots, m$, of ${\bar M}_c$ such that the following holds for the maximum likelihood estimator ${\hat{\bar u}}^{{\bar a}}_{mle}$ of 
${\bar u}^{{\bar a}}$ without bias-correction:
\begin{align}
\label{eq:36}
E[K\nu({\hat{\bar u}}^{{\bar a}}_{mle} - {\bar u}^{\bar a})
({\hat{\bar u}}^{{\bar b}}_{mle} -{\bar u}^{\bar b})] = 
g^{{\bar a}{\bar b}} + O(K^{-2}).
\end{align}
When $M_c$ itself is a f.r.m.\ exponential family, (\ref{eq:36}) holds by 
(\ref{eq:21}) given in Theorem 4.1. 
When $M_c$ is a dual quadric hypersurface, 
(\ref{eq:36}) holds by (\ref{eq:29}) given in Theorem 5.3.  
\end{theorem}
\

\begin{proof}
Since $H_{\kappa \lambda a}^{(-1)} = 0$ holds for the maximum likelihood estimator (m.l.e.), from (\ref{eq:34}) and Theorem 5.4, we have for the bias of the m.l.e.\
\begin{align*}
b_{mle}^{\bar a} = - \ ^{'}\!\Gamma_{\alpha \beta}^{(-1){\bar a}}
g^{\alpha \beta} = 
- \ ^{'}\!\Gamma_{{\bar b}{\bar c}}^{(-1){\bar a}}g^{{\bar b}{\bar c}} 
- \ ^{'}\!\Gamma_{\kappa \lambda}^{(-1){\bar a}}g^{\kappa \lambda} = 
- \ ^{'}\!\Gamma_{{\bar b}{\bar c}}^{(-1){\bar a}}g^{{\bar b}{\bar c}} 
- H_{\kappa \lambda}^{(-1){\bar a}}g^{\kappa \lambda} = 0,
\end{align*}
and the expression (\ref{eq:36}) is derived.

When $M_c$ itself is a f.r.m.\ exponential family with expectation parameter $u$, the partial differential equations for $s(u) = \log \nu(u)$ and for ${\bar u} = {\bar u}(u)$ are
\begin{align*}
\partial_a s_b - s_as_b = 0,\quad 
\partial_bC_a^{\bar a} - s_bC_a^{\bar a} - s_aC_b^{\bar a} = 0,\ \ 
C_a^{\bar a} = \frac{\partial {\bar u}^{\bar a}}{\partial u^a},
\end{align*}
as noted in Theorem 4.1. 
When $M_c$ is a dual quadric hypersurface, the partial differential equations for $s(u) = \log \nu(u)$ and for ${\bar u} = {\bar u}(u)$ are
\begin{align*}
\partial_a s_b - \Gamma_{ab}^{(-1)c}s_c - s_as_b = k_0l_0g_{ab},\quad
\partial_bC_a^{\bar a} - 
(\Gamma_{ba}^{(-1)c} + \delta_a^cs_b + \delta_b^cs_a)
C_c^{\bar a} = 0,
\end{align*}
as noted in Theorem 5.3.

This completes the proof of the theorem.
\end{proof}

\section{Examples}

\subsection{von Mises-Fisher model}

This is an $(m+1, m)$-curved exponential family, of which density functions with respect to the invariant measure on the $m$-dimensional unit sphere under rotational transformations are given by 
(cf.\ Barndorff-Nielsen et al (1989), p.76)
\begin{align*}
&f_c(x, 1, u) = \exp\{\theta(u)\cdot x - \psi(\theta(u))\},\quad 
\theta\cdot x = \theta^1x_1 + \theta^2x_2 + \cdots + \theta^{m+1}x_{m+1},\\ 
&\theta = r\xi = (r\xi^i),\ 
\xi \in S^m = \{\xi \in {\mathbb R}^{m+1}\ |\ 
\xi\cdot \xi = 1\},\quad 
x = (x_i) \in S^m,\\
&\psi(\theta) = -\log a_m(r),\quad
1/a_m(r) = (2\pi)^{(m+1)/2}r^{(1-m)/2}I_{(m-1)/2}(r),\quad r > 0,
\end{align*}
where $I_{(m-1)/2}(r)$ is the modified Bessel function of the first kind and of order $(m-1)/2$. We assume that the concentration parameter $r$ is assumed to be a given positive constant. The parametric representations $\theta = \theta(u)$ and $\eta = \eta(u)$ are given by\\

\begin{tabular}{ll}
$\theta^1(u) = r\cos u^1$ & $\eta_1(u) = r^\dagger \cos u^1$ \\
$\theta^2(u) = r\sin u^1 \cos u^2$ & 
$\eta_2(u) = r^\dagger \sin u^1 \cos u^2$ \\
$\theta^3(u) = r\sin u^1 \sin u^2 \cos u^3$ & 
$\eta_3(u) = r^\dagger \sin u^1 \sin u^2 \cos u^3$ \\
$\qquad \qquad \cdots$ & $\qquad \qquad \cdots$ \\
$\theta^{m+1}(u) = r\sin u^1 \sin u^2 \cdots \sin u^{m-1} \sin u^m$ & 
$\eta_{m+1}(u) = r^\dagger \sin u^1 \sin u^2 \cdots \sin u^{m-1} 
\sin u^m$,  
\end{tabular}
\

\noindent
where $0 \le u^1, \dots, u^{m-1} \le \pi,\ 0 \le u^m < 2\pi$ and
$r^\dagger = -d \log a_m(r)/dr = I_{(m+1)/2}(r)/I_{(m-1)/2}(r)$. 
Note that $E[x] = r^\dagger \xi$, 
and $r^\dagger$ is a strictly increasing function of $r$ which maps 
$(0, \infty)$ onto $(0, 1)$.

From these representations the tangent vectors $B_a^i(u)$ and $B_{ai}(u)$ can be calculated, and then the unit normal vectors 
$B_\kappa^i(u)$ and $B_{\kappa i}(u)\ (\kappa = m+1)$ are derived from the relations 
$B_\kappa^i(u)B_{ai}(u) = 0$ and $B_{\kappa i}(u)B_a^i(u) = 0$ 
as follows. \\

\begin{tabular}{ll}
$B_\kappa^1(u) = \cos u^1$ & $B_{\kappa 1}(u) = \cos u^1$ \\
$B_\kappa^2(u) = \sin u^1 \cos u^2$ & 
$B_{\kappa 2}(u) = \sin u^1 \cos u^2$ \\
$B_\kappa^3(u) = \sin u^1 \sin u^2 \cos u^3$ & 
$B_{\kappa 3}(u) = \sin u^1 \sin u^2 \cos u^3$ \\
$\qquad \qquad \cdots$ & $\qquad \qquad \cdots$ \\
$B_\kappa^{m+1}(u) = \sin u^1 \sin u^2 \cdots \sin u^{m-1} \sin u^m$ & 
$B_{\kappa\ {m+1}}(u) = \sin u^1 \sin u^2 \cdots \sin u^{m-1} 
\sin u^m$.  
\end{tabular}
\

\noindent
The above expressions show that
\begin{align*}
B_\kappa^i(u) = \frac{1}{r}\theta^i(u),\quad 
B_{\kappa i}(u) = \frac{1}{r^\dagger} \eta_i(u),
\end{align*}
and so this model is a dual quadric hypersurface. 
From Theorem 5.2 we also see that this model is ES conjugate 
$m$($e$)-symmetric and totally $e$-umblic with constant $m$($e$)-curvature 
$1/(rr^\dagger) > 0$. 
The related geometrical quantities are given below. 
\begin{align*}
&g_{ab}(u) = \delta_{ab}rr^\dagger \prod_{c=1}^a\sin^2 u^{c-1},\quad 
\sin^2 u^0 = 1,\\
&H_{ab\kappa}^{(1)}(u) = -\frac{1}{r^\dagger}g_{ab}(u),\quad
H_{ab\kappa}^{(-1)}(u) = -\frac{1}{r}g_{ab}(u) = 
\frac{r^\dagger}{r}H_{ab\kappa}^{(1)}(u),\\
&R_{abba}^{(\pm1)}(u) = \frac{1}{rr^\dagger}g_{aa}(u)g_{bb}(u),\quad 
a \neq b,\quad R_{ab}^{(\pm1)}(u) = \frac{m-1}{rr^\dagger}g_{ab}(u).
\end{align*}

From Theorem 3.2 (i) this $M_c$ is conformally $m$($e$)-flat, so that 
there exist a gauge function $\nu(u) > 0$ and a $(-1)$-affine coordinate system ${\bar u} = ({\bar u}^{{\bar a}})$ such that  
${\bar \Gamma}_{{\bar a}{\bar b}{\bar c}}^{(-1)}({\bar u}) = 0,\ 
\forall {\bar u}\in {\bar M}_c$. 
As given by (\ref{eq:29}), the partial differential equation for $s(u) = \log \nu(u)$ is
\begin{align*}
\partial_as_b - \Gamma_{ab}^{(-1)c}s_c - s_as_b = 
\frac{1}{rr^\dagger}g_{ab},
\end{align*}
of which one solution is
\begin{align*}
\nu(u) = \frac{1}{\prod_{a=1}^m |\sin u^a|},
\end{align*}
and then 
${\bar u}^{\bar a} = \nu(u)D^{{\bar a}i}\eta_i(u)$.

\subsection{Hyperboloid model}

This is an $(m+1, m)$-curved exponential family, of which density functions with respect to the invariant measure on the $m$-dimensional unit hyperboloid under hyperbolic transformations are given by 
(cf.\ Barndorff-Nielsen et al. (1989), p.104)
\begin{align*}
&f_c(x, 1, u) = \exp\{\theta(u)\cdot x - \psi(\theta(u))\},\quad 
\theta^1 = -r\xi^1,\ \ \theta^i = r\xi^i,\ i = 2, \dots, m+1,\\
&\xi = (\xi^i) \in H^m = \{\xi \in {\mathbb R}^{m+1}\ |\ 
\xi*\xi = 1,\ \xi^1 > 0 \},\quad 
x = (x_i) \in H^m,\\
&\xi*\xi = (\xi^1)^2 - (\xi^2)^2 - \cdots - (\xi^{m+1})^2,\\ 
&\psi(\theta) = -\log a_m(r),\quad
1/a_m(r) = 2(2\pi)^{(m-1)/2}r^{(1-m)/2}K_{(m-1)/2}(r),\quad r > 0,
\end{align*}
where $K_{(m-1)/2}(r)$ is the modified Bessel function of the third kind and of order $(m-1)/2$. We assume that the concentration parameter $r$ is assumed to be a given positive constant. The parametric representations $\theta = \theta(u)$ and $\eta = \eta(u)$ are given by\\

\begin{tabular}{ll}
$\theta^1(u) = -r\cosh u^1$ & $\eta_1(u) = r^\dagger \cosh u^1$ \\
$\theta^2(u) = r\sinh u^1 \cos u^2$ & 
$\eta_2(u) = r^\dagger \sinh u^1 \cos u^2$ \\
$\theta^3(u) = r\sinh u^1 \sin u^2 \cos u^3$ & 
$\eta_3(u) = r^\dagger \sinh u^1 \sin u^2 \cos u^3$ \\
$\qquad \qquad \cdots$ & $\qquad \qquad \cdots$ \\
$\theta^{m+1}(u) = r\sinh u^1 \sin u^2 \cdots \sin u^{m-1} \sin u^m$ & 
$\eta_{m+1}(u) = r^\dagger \sinh u^1 \sin u^2 \cdots \sin u^{m-1} 
\sin u^m$,  
\end{tabular}
\

\noindent
where $u^1 \in {\mathbb R},\ 0 \le u^2, \dots, u^{m-1} \le \pi,\ 
0 \le u^m < 2\pi$ and 
$r^\dagger = d \log a_m(r)/dr = K_{(m+1)/2}(r)/K_{(m-1)/2}(r)$. 
Note that $E[x] = r^\dagger \xi$, 
and $r^\dagger$ is a strictly decreasing function of $r$ which maps 
$(0, \infty)$ onto $(1, \infty)$. 

From these representations the tangent vectors $B_a^i(u)$ and $B_{ai}(u)$ can be calculated, and then the unit normal vectors 
$B_\kappa^i(u)$ and $B_{\kappa i}(u)\ (\kappa = m+1)$ are derived from the relations 
$B_\kappa^i(u)B_{ai}(u) = 0$ and $B_{\kappa i}(u)B_a^i(u) = 0$ 
as follows. \\

\begin{tabular}{ll}
$B_\kappa^1(u) = \cosh u^1$ & $B_{\kappa 1}(u) = \cosh u^1$ \\
$B_\kappa^2(u) = -\sinh u^1 \cos u^2$ & 
$B_{\kappa 2}(u) = \sinh u^1 \cos u^2$ \\
$B_\kappa^3(u) = -\sinh u^1 \sin u^2 \cos u^3$ & 
$B_{\kappa 3}(u) = \sinh u^1 \sin u^2 \cos u^3$ \\
$\qquad \qquad \cdots$ & $\qquad \qquad \cdots$ \\
$B_\kappa^{m+1}(u) = -\sinh u^1 \sin u^2 \cdots \sin u^{m-1} \sin u^m$ & 
$B_{\kappa\ {m+1}}(u) = \sinh u^1 \sin u^2 \cdots \sin u^{m-1} 
\sin u^m$.  
\end{tabular}
\

\noindent
The above expressions show that
\begin{align*}
B_\kappa^i(u) = -\frac{1}{r}\theta^i(u),\quad 
B_{\kappa i}(u) = \frac{1}{r^\dagger} \eta_i(u),
\end{align*}
and again this model is a dual quadric hypersurface. 
From Theorem 5.2 we also see that this model is ES conjugate 
$m$($e$)-symmetric and totally $e$-umblic with constant $m$($e$)-curvature 
$-1/(rr^\dagger) < 0$. 
The related geometrical quantities are given below. 
\begin{align*}
&g_{11}(u) = rr^\dagger,\quad 
g_{ab}(u) = \delta_{ab}rr^\dagger \sinh^2 u^1 \prod_{c=2}^a\sin^2 u^{c-1},
\ \ \sin^2 u^1 = 1,\ \ a = 2, \dots, m,\\
&H_{ab\kappa}^{(1)}(u) = -\frac{1}{r^\dagger}g_{ab}(u),\quad
H_{ab\kappa}^{(-1)}(u) = \frac{1}{r}g_{ab}(u) = 
-\frac{r^\dagger}{r}H_{ab\kappa}^{(1)}(u),\\
&R_{abba}^{(\pm1)}(u) = -\frac{1}{rr^\dagger}g_{aa}(u)g_{bb}(u),\quad 
a \neq b,\quad R_{ab}^{(\pm1)}(u) = -\frac{m-1}{rr^\dagger}g_{ab}(u).
\end{align*}

From Theorem 3.2 (i) this $M_c$ is conformally $m$($e$)-flat, so that 
there exist a gauge function $\nu(u) > 0$ and a $(-1)$-affine coordinate system ${\bar u} = ({\bar u}^{{\bar a}})$ such that  
${\bar \Gamma}_{{\bar a}{\bar b}{\bar c}}^{(-1)}({\bar u}) = 0,\ 
\forall {\bar u}\in {\bar M}_c$. 
As given by (\ref{eq:29}), the partial differential equation for $s(u) = \log \nu(u)$ is
\begin{align*}
\partial_as_b - \Gamma_{ab}^{(-1)c}s_c - s_as_b = 
-\frac{1}{rr^\dagger}g_{ab},
\end{align*}
of which one solution is
\begin{align*}
\nu(u) = \frac{1}{|\sinh u^1|\prod_{a=2}^m |\sin u^a|},
\end{align*}
and then 
${\bar u}^{\bar a} = \nu(u)D^{{\bar a}i}\eta_i(u)$.

\subsection{Numerical results}

We examine our theoretical results numerically by using the von Mises-Fisher and the hyperboloid models. 
We take 10 kinds of number $N$ (nonsequential case) and $K$ (sequential case)  of observations, and for each $N$ or $K$, we generate $500$ random simulated data. 
Then the empirical means of covariances  
$E_{emp}[({\hat u}^{*a}_{mle} - u_0^a)({\hat u}^{*b}_{mle} - u_0^b)]$ (nonsequential case) and 
$E_{emp}[({\hat {\bar u}}^{\bar a}_{mle} - {\bar u}_0^{\bar a})
({\hat {\bar u}}^{\bar b}_{mle} - {\bar u}_0^{\bar b})]$ 
(sequential case) of the m.l.e.\ over this 500 sample size are used for evaluation, where $u_0^a$ and ${\bar u}_0^{\bar a}$ denote the true values of $u^a$ and ${\bar u}^{\bar a}$. The stopping times $\tau$ for the sequential estimations are determined by 
(see Okamoto, Amari and Takeuchi (1991))
\begin{align*}
&\tau = \inf \Big\{t\ |\ -\frac{1}{m}\partial_a\partial_b 
l(x, t, {\hat u}_{mle})g^{ab}({\hat u}_{mle}) \ge K\nu({\hat u}_{mle}) + 
c \Big\},\\ 
&c = -\frac{1}{2}\Big(\frac{m}{rr^\dagger} - \frac{1}{r^{\dagger 2}} \Big)\ :\ \textrm{von Mises-Fisher}\quad 
c = -\frac{1}{2}\Big(-\frac{m}{rr^\dagger} - \frac{1}{r^{\dagger 2}} \Big)\ :\ \textrm{hyperboloid}.
\end{align*}

As for the von Mises-Fisher model, numerical results are based on the following set of values
\begin{align*}
&m = 2,\quad r = 0.25,\quad 
(u_0^1, u_0^2) = (\pi/6, \pi/3),\\
&D^{{\bar a}i} = \delta^{{\bar a}i},\ \textrm{i.e.,}\ 
{\bar u}^1 = \nu(u)\eta_1(u),\quad 
{\bar u}^2 = \nu(u)\eta_2(u),\quad 
\nu(u) = 1/(|\sin u^1|| \sin u^2|),
\end{align*}
and for the hyperboloid model, numerical results are based on the following set of values
\begin{align*}
&m = 2,\quad r = 0.1,\quad 
(u_0^1, u_0^2) = (0.1, \pi/3),\\
&D^{{\bar a}i} = \delta^{{\bar a}i}/100,\ \textrm{i.e.,}\ 
{\bar u}^1 = \nu(u)\eta_1(u)/100,\quad 
{\bar u}^2 = \nu(u)\eta_2(u)/100,\quad 
\nu(u) = 1/(|\sinh u^1|| \sin u^2|).
\end{align*}

Figures 3-8 show the von Mises-Fisher model, and Figures 9-14 show the hyperboloid model. 
The notations in the figures indicate the following quantities.
\begin{align*}
&\textrm{$\cdot$ nonsequential case}\\
&\quad OCOVab = NE_{emp}[({\hat u}^{*a}_{mle} - u_0^a)({\hat u}^{*b}_{mle} - u_0^b)],\quad 
a, b = 1, 2\\
&\quad OCRBab = g^{ab}(u_0),\\
&\quad OALBab = g^{ab}(u_0) + \frac{1}{N}\bigg\{\frac{1}{2} (\Gamma_{M_c}^{(-1)}(u_0))^{2ab} + 
 (H_{M_c}^{(1)}(u_0))^{2ab}  \bigg\}\\
&\textrm{$\cdot$ sequential case}\\
&\quad CCOV{\bar a}{\bar b} = 
E_{emp}(\tau)E_{emp}[({\hat {\bar u}}^{\bar a}_{mle} - {\bar u}_0^{\bar a})
({\hat {\bar u}}^{\bar b}_{mle} - {\bar u}_0^{\bar b})],\quad 
{\bar a}, {\bar b} = 1, 2\\
&\quad CCRB{\bar a}{\bar b} = g^{{\bar a}{\bar b}}({\bar u}_0),\\
&\quad MST = E_{emp}(\tau)\ :\ \textrm{empirical mean of $\tau$}\\
&\quad SDST = \sqrt{V_{emp}(\tau)}\ :\ \textrm{empirical standard deviation of $\tau$}.
\end{align*}

We see that in the nonsequential case $OCOVab$ approach to the asymptotic lower bound $OALBab$ exhibiting the differential geometrical loss  
$OALBab - OCRBab$, and in the sequential case $CCOV{\bar a}{\bar b}$ nearly attain the Cram\'er-Rao lower bound $CCRB{\bar a}{\bar b}$ as if the model were a f.r.m.\ exponential family. 
Figures 8, 14 confirm that the assumptions $MST= O(K), SDST = O(\sqrt{K})$ 
are satisfied in each model.

\clearpage

\begin{figure}[htbp]
\begin{minipage}{.5\linewidth}
\begin{center}
\includegraphics[width=8cm,height=7cm]{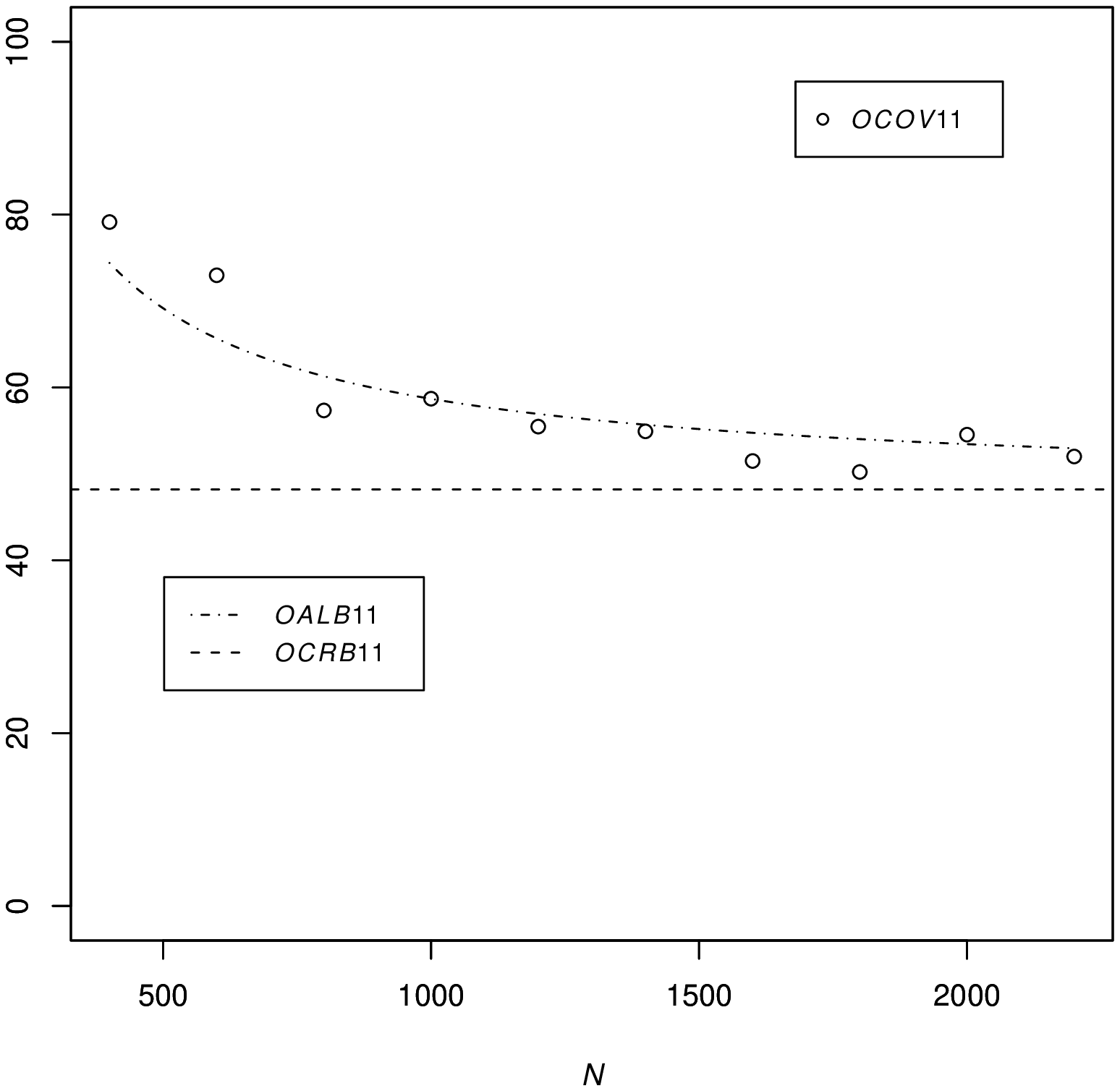}
\vspace*{-11mm}   
\caption{$OCOV11$\ \ von Mises-Fisher}
\label{fig:1-1}
\end{center}
\end{minipage}
\begin{minipage}{.5\linewidth}
\begin{center}
\includegraphics[width=8cm,height=7cm]{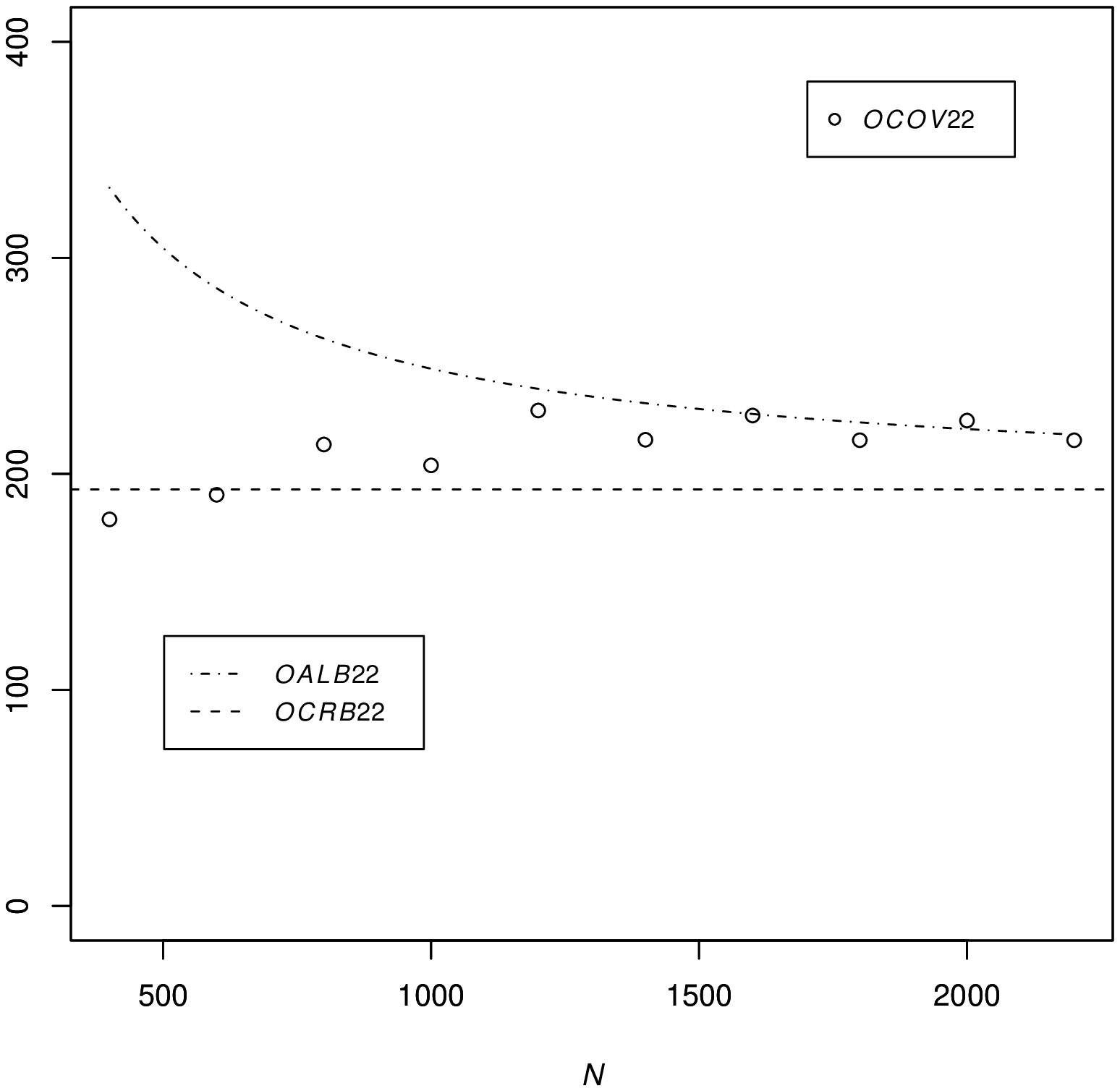}
\vspace*{-11mm}   
\caption{$OCOV22$\ \ von Mises-Fisher}
\label{fig:1-1}
\end{center}
\end{minipage}
\begin{minipage}{.5\linewidth}
\begin{center}
\includegraphics[width=8cm,height=7cm]{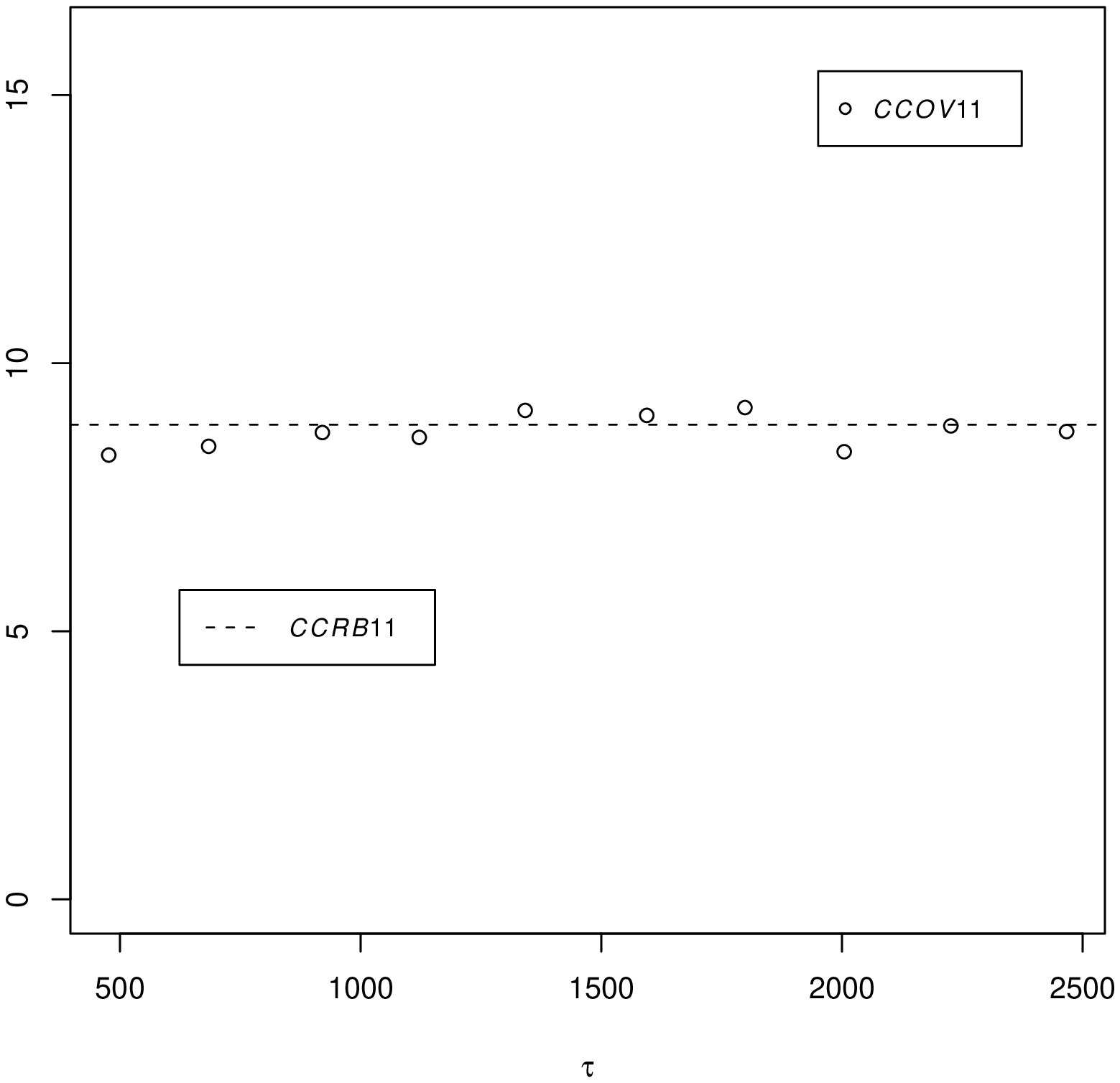}
\vspace*{-11mm}   
\caption{$CCOV11$\ \ von Mises-Fisher}
\label{fig:1-2}
\end{center}
\end{minipage}
\begin{minipage}{.5\linewidth}
\begin{center}
\includegraphics[width=8cm,height=7cm]{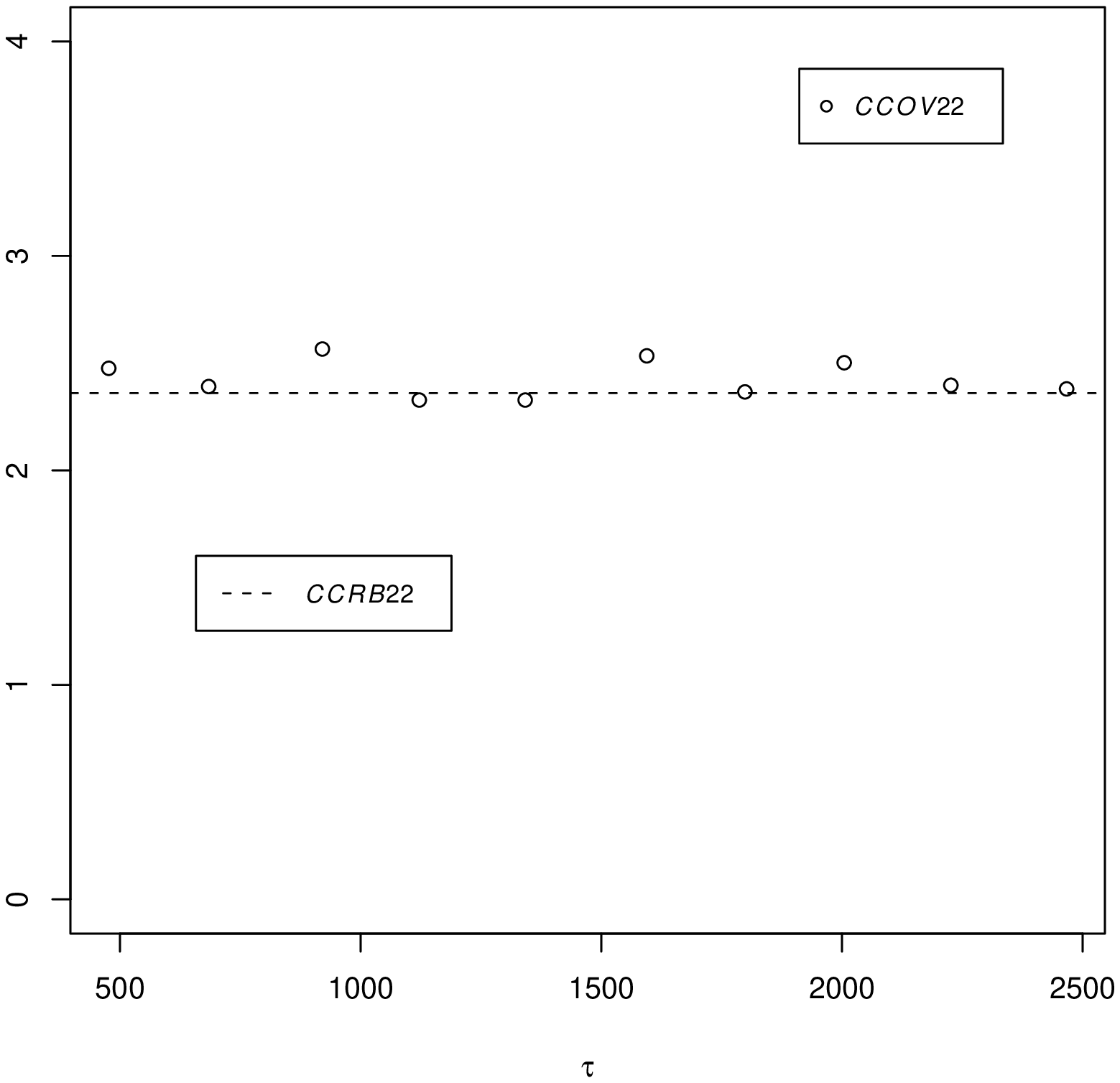}
\vspace*{-11mm}   
\caption{$CCOV22$\ \ von Mises-Fisher}
\label{fig:1-3}
\end{center}
\end{minipage}
\begin{minipage}{.5\linewidth}
\begin{center}
\includegraphics[width=8cm,height=7cm]{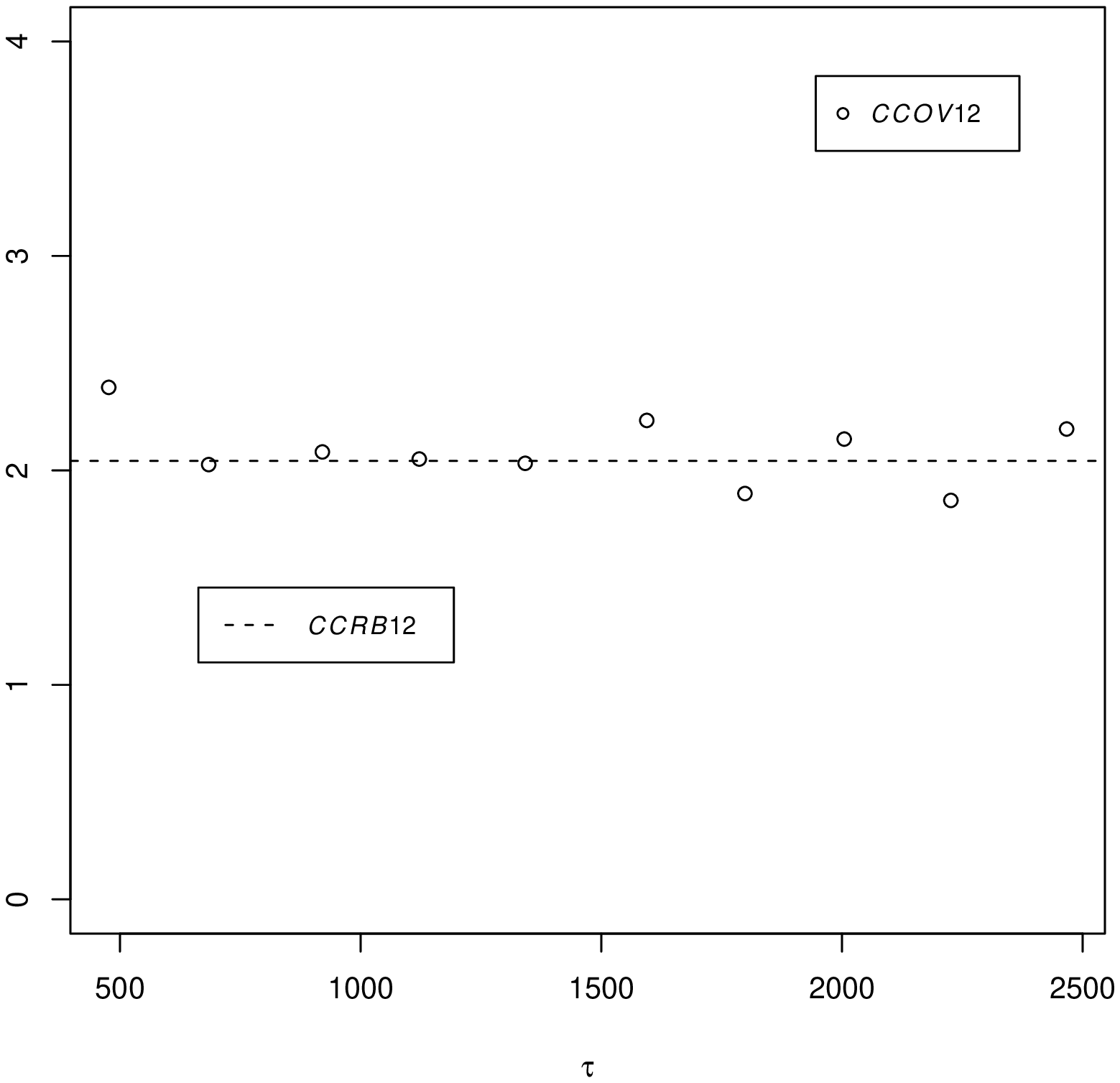}
\vspace*{-11mm}   
\caption{$CCOV12$\ \ von Mises-Fisher}
\label{fig:1-4}
\end{center}
\end{minipage}
\begin{minipage}{.5\linewidth}
\begin{center}
\includegraphics[width=8cm,height=7cm]{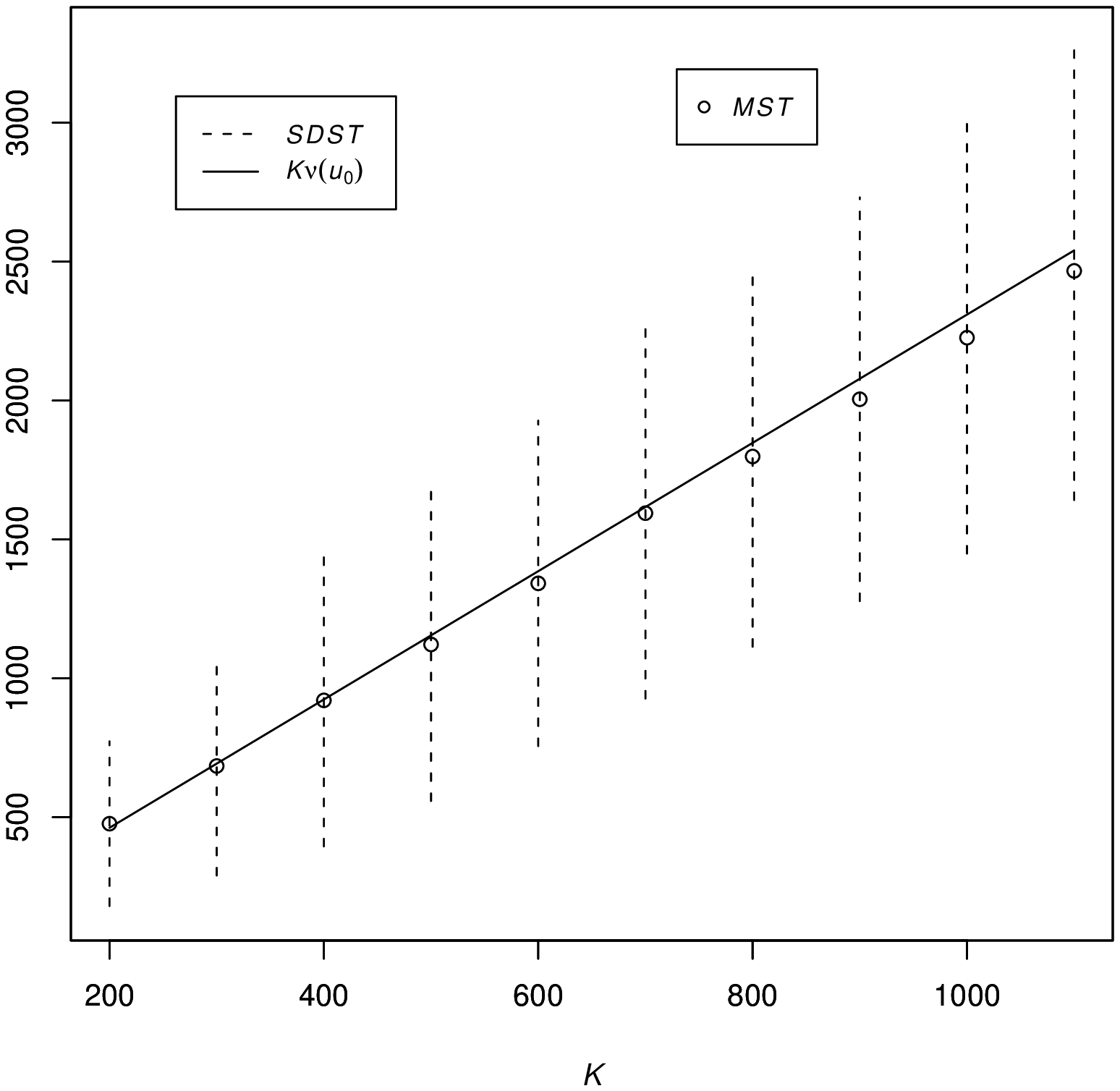}
\vspace*{-11mm}   
\caption{$MST\ SDST$\ \ von Mises-Fisher}
\label{fig:1-5}
\end{center}
\end{minipage}
\end{figure}

\clearpage

\begin{figure}[htbp]
\begin{minipage}{.5\linewidth}
\begin{center}
\includegraphics[width=8cm,height=7cm]{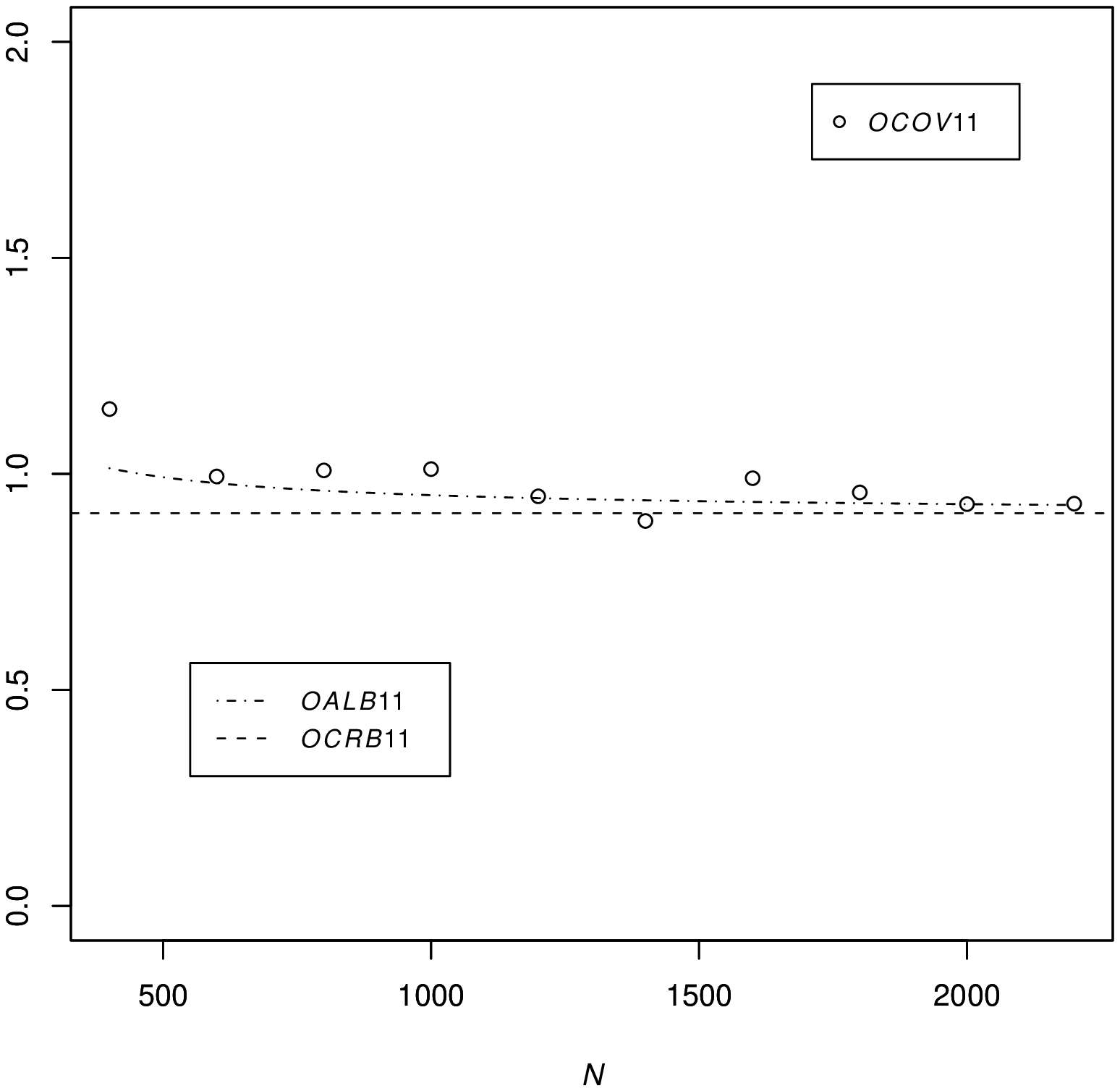}
\vspace*{-11mm}   
\caption{$OCOV11$\ \ hyperboloid}
\label{fig:1-1}
\end{center}
\end{minipage}
\begin{minipage}{.5\linewidth}
\begin{center}
\includegraphics[width=8cm,height=7cm]{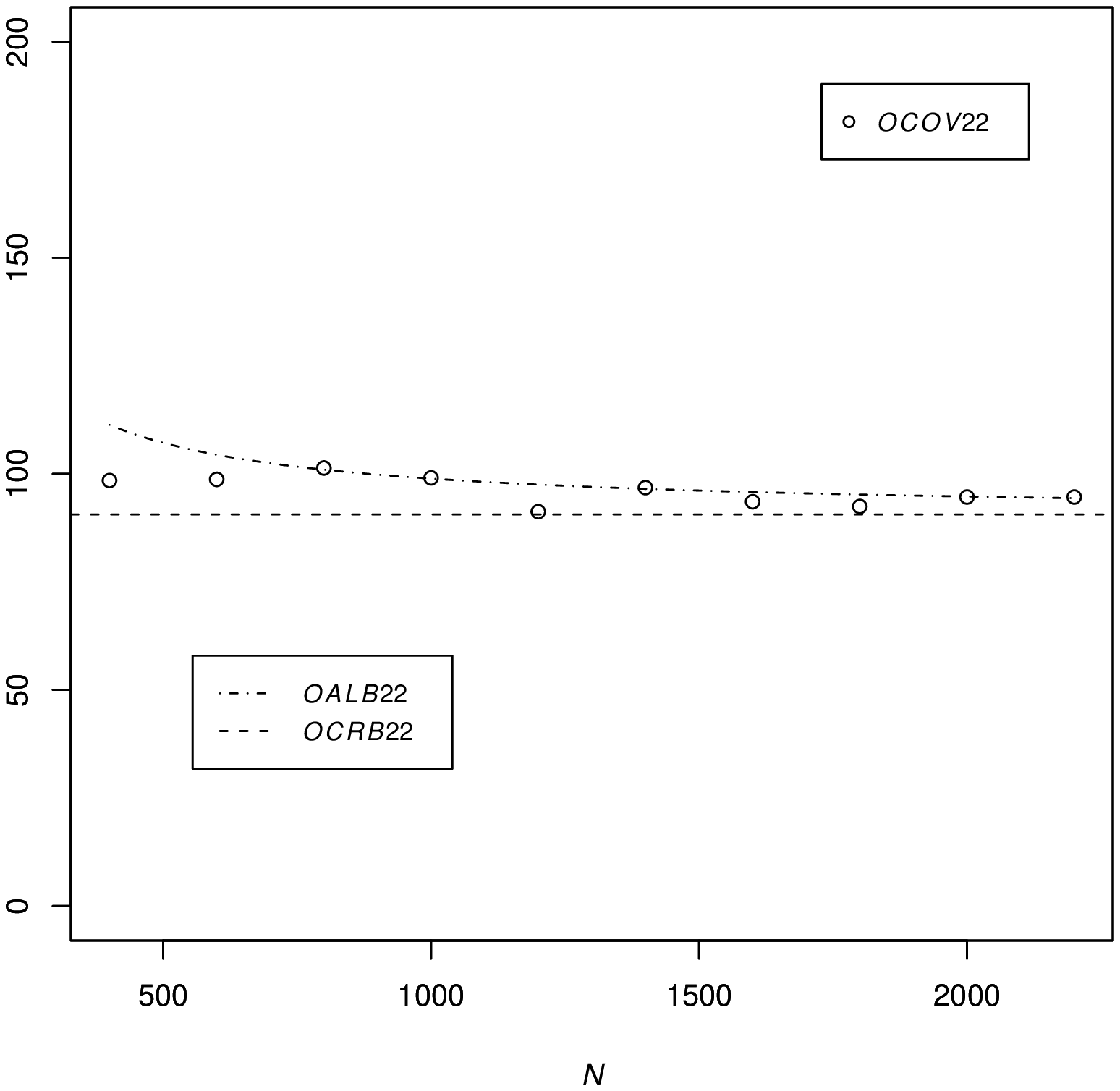}
\vspace*{-11mm}   
\caption{$OCOV22$\ \ hyperboloid}
\label{fig:1-1}
\end{center}
\end{minipage}
\begin{minipage}{.5\linewidth}
\begin{center}
\includegraphics[width=8cm,height=7cm]{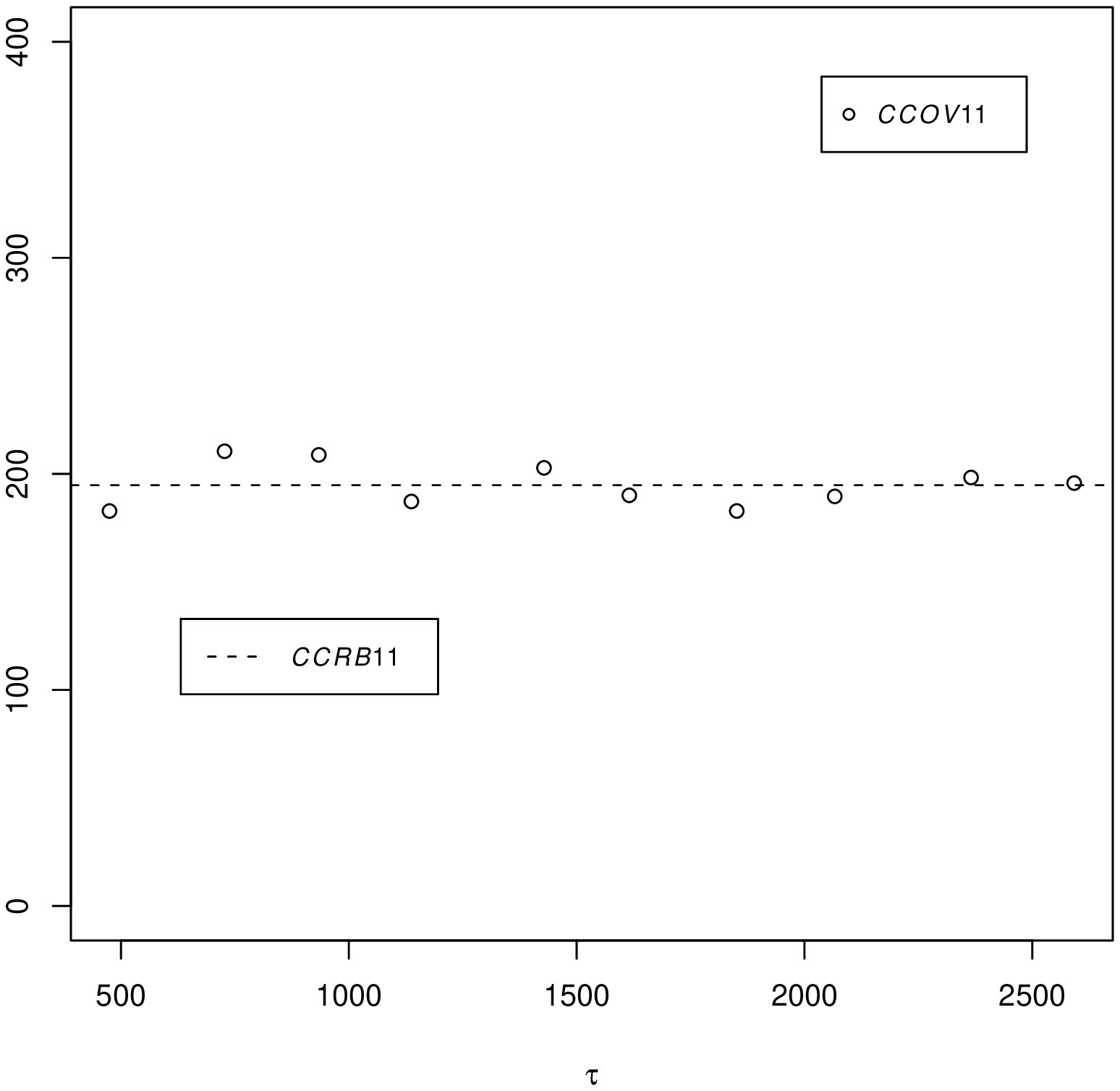}
\vspace*{-11mm}   
\caption{$CCOV11$\ \ hyperboloid}
\label{fig:1-2}
\end{center}
\end{minipage}
\begin{minipage}{.5\linewidth}
\begin{center}
\includegraphics[width=8cm,height=7cm]{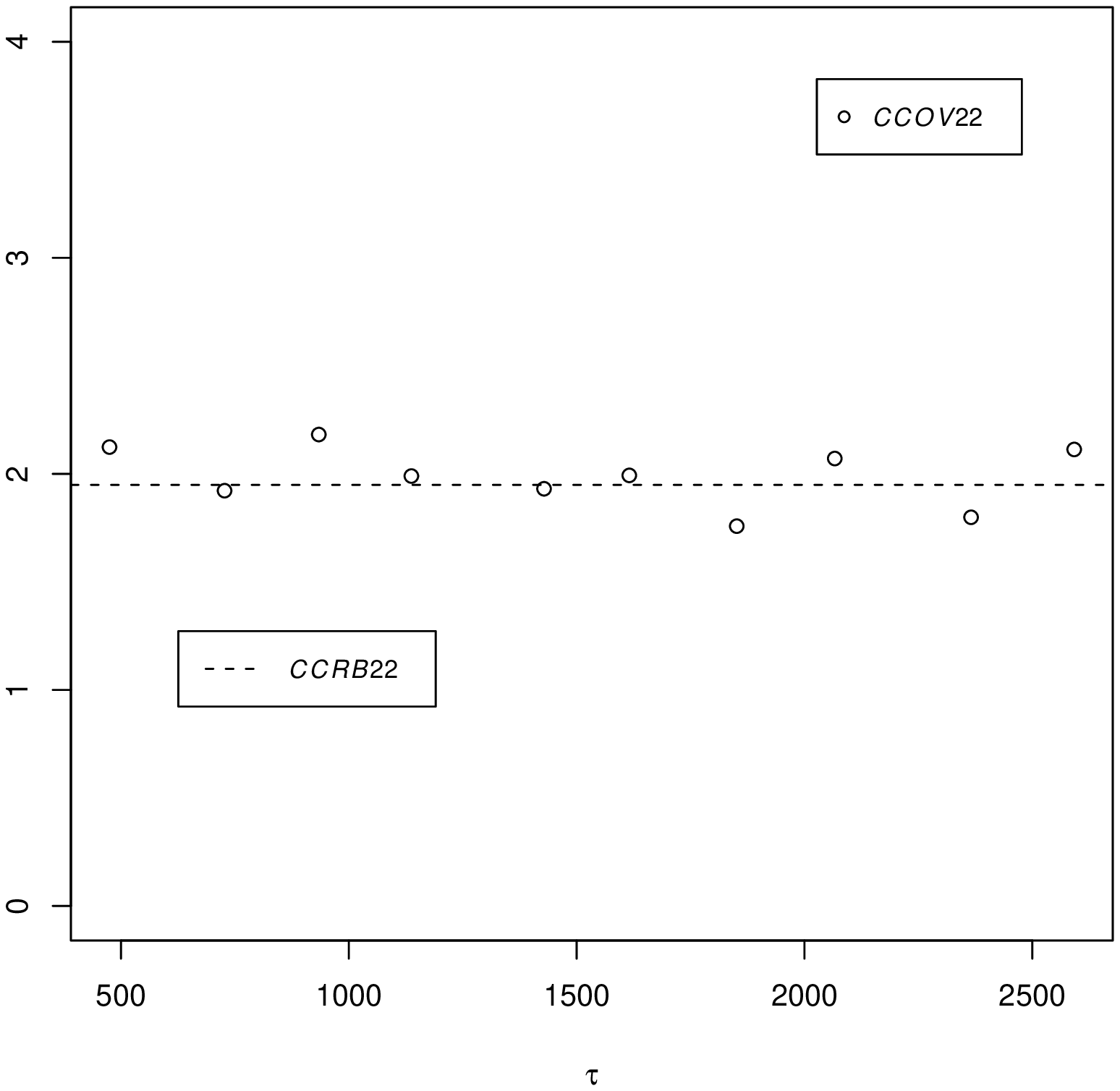}
\vspace*{-11mm}   
\caption{$CCOV22$\ \ hyperboloid}
\label{fig:1-3}
\end{center}
\end{minipage}
\begin{minipage}{.5\linewidth}
\begin{center}
\includegraphics[width=8cm,height=7cm]{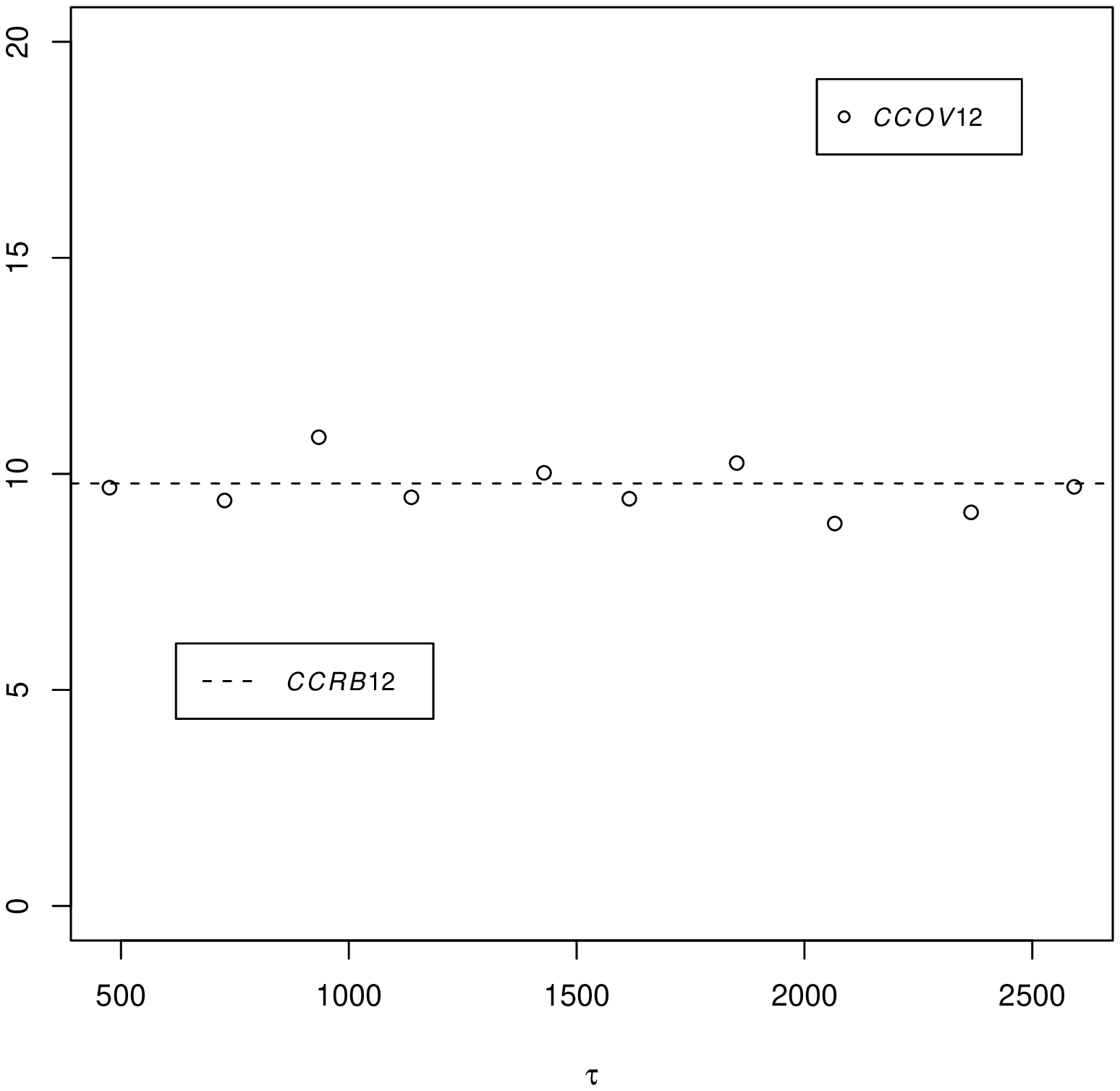}
\vspace*{-11mm}   
\caption{$CCOV12$\ \ hyperboloid}
\label{fig:1-4}
\end{center}
\end{minipage}
\begin{minipage}{.5\linewidth}
\begin{center}
\includegraphics[width=8cm,height=7cm]{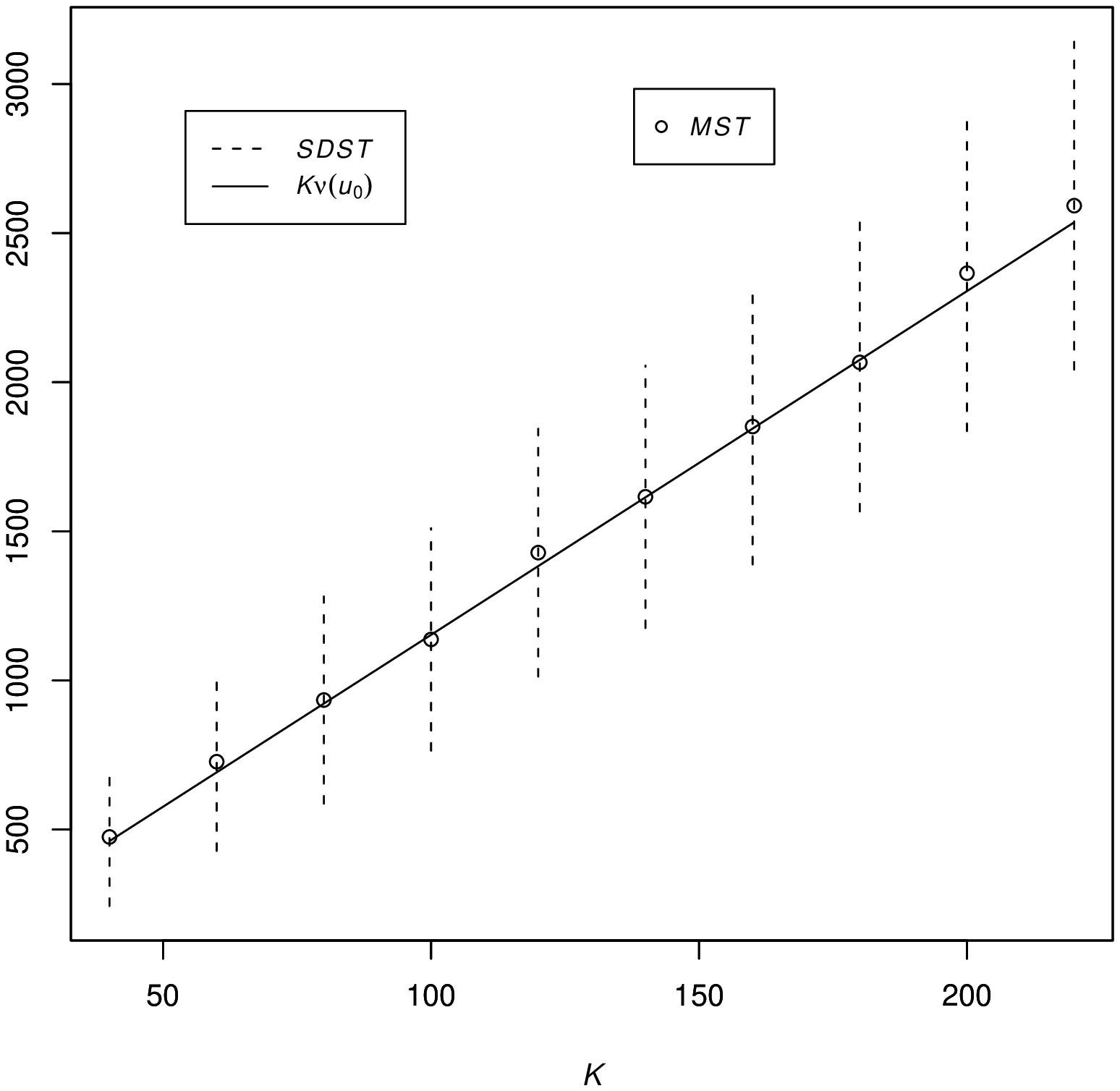}
\vspace*{-11mm}   
\caption{$MST\ SDST$\ \ hyperboloid}
\label{fig:1-5}
\end{center}
\end{minipage}
\end{figure}

\clearpage

\section{Discussion}

We have analyzed sequential estimation procedures in terms of the conformal geometry of statistical manifolds. We have also constructed a  concrete procedure for the covariance mininization in a multidimensional curved exponential family $M_c$. The method is divided into two separate stages: one is to choose a stopping rule which is effective for reducing the $1$-ES curvature $H^{(1)}_{M_c}$ and the other is to choose a gauge function $\nu(u)$ on $M_c$ effective for reducing the $(-1)$-connection 
$\Gamma^{(-1)}_{M_c}$. Another typical choice of $\nu(u)$ is the one effective for the covariance stabilization, as suggested in Okamoto, Amari and Takeuchi (1991). These choices contradict each other in general multidimensional cases, and this fact reflects the difference between the ordinary Riemannian geometry and the mutually dual geometry as exhibited in several geometrical notions introduced in this paper.

The present method is also applicable to investigating sequential testing procedures. The geometrical theory of higher-order asymptotics of testing hypothesis in nonsequential case was developed by Kumon and Amari (1983) and Amari (1985). The main results are summarized as follows.

The power function $P_T(t)$ of a test $T$ is expanded as
\begin{align*}
P_T(t) = P_{T1}(t) + P_{T2}(t)/\sqrt{N} + P_{T3}(t)/N + O(N^{-3/2}),
\end{align*}
where $N$ denotes the number of observations, and $t/\sqrt{N}$ indicates  the geodesic distance between the null hypothesis and the point in the alternative hypothesis.\\
\quad (i)\ The first-order power function $P_{T1}(t)$ and the second-order power function $P_{T2}(t)$ are maximized uniformly in $t$ if and only if the ancillary family (boundaries of the critical region) associated with a test $T$ is asymptotically an orthogonal family. \\
\quad (ii)\ The third-order power loss function
$\Delta P_{T3}(t) = \sup_T P_{T3}(t) - P_{T3}(t)$ is expressed as the weighted sum of two kinds of the square of the $1$-ES curvatures $H^{(1)}_{M_c}$, the square of the $(-1)$-ES mixture curvature 
$H^{(-1)}_T$ of the associated ancillary family, and also the square of the $(-1)$-mixture connection $\Gamma^{(-1)}_{M_c}$ (when there are unknown nuisance parameters).

Based on these nonsequential results, we can utilize the conformal geometry to the analysis and the construction of most powerful sequential tests. Specifically when a statistical manifold is a f.r.m. exponential family or a dual quadric hypersurface, it is expected that one can design sequential tests without any power loss. This is a subject which will be treated in a future work.

\end{document}